\documentclass[11pt,reqno]{amsart}
\usepackage[T1]{fontenc}
\usepackage{amsmath,amssymb,amsthm,mathrsfs}
\usepackage[margin=1in]{geometry}
\usepackage[hidelinks]{hyperref}
\numberwithin{equation}{section}
\newtheorem{theorem}{Theorem}[section]
\newtheorem{proposition}[theorem]{Proposition}
\newtheorem{lemma}[theorem]{Lemma}
\newtheorem{corollary}[theorem]{Corollary}
\theoremstyle{remark}
\newtheorem{remark}[theorem]{Remark}
\usepackage[colorinlistoftodos]{todonotes}
\usepackage{etoolbox}

\makeatletter 
\@mparswitchfalse
\makeatother
\normalmarginpar 

\makeatletter
\patchcmd{\@settitle}
{\uppercasenonmath\@title}
{}{}{}
\patchcmd{\@settitle}
{\bfseries}
{\fontsize{16}{20}\selectfont\bfseries}
{}{}

\patchcmd{\@setauthors}
{\MakeUppercase{\authors}}
{\authors}{}{}
\patchcmd{\@setauthors}
{\footnotesize}
{\fontsize{11}{15}\selectfont}
{}{}

\renewcommand{\author@andify}{%
	\nxandlist
	{\unskip\quad}
	{\unskip\quad}
	{\unskip\quad}%
}

\patchcmd{\@setaddresses}
{\indent}
{\noindent}
{}{}

\patchcmd{\@setaddresses}
{\indent\emailaddrname}
{\noindent\emailaddrname}
{}{}
\makeatother

\newcommand{\R}{\mathbb R}
\newcommand{\C}{\mathbb C}
\newcommand{\HH}{\mathbb H}
\newcommand{\CP}{\mathbb{CP}}
\newcommand{\HP}{\mathbb{HP}}
\newcommand{\eps}{\varepsilon}

\newcommand{\Zcal}{\mathcal Z}
\newcommand{\Kcal}{\mathcal K}
\newcommand{\Acal}{\mathcal A}
\newcommand{\Pcal}{\mathcal P}
\newcommand{\dd}{\,d\mathrm{vol}_g}
\newcommand{\rp}[2]{\langle #1,#2\rangle}
\newcommand{\ga}{\gamma}

\DeclareMathOperator{\Crit}{Crit}
\DeclareMathOperator{\ind}{ind}
\DeclareMathOperator{\Cat}{Cat}
\DeclareMathOperator{\Ran}{Ran}
\DeclareMathOperator{\Scal}{Scal}
\DeclareMathOperator{\Vol}{Vol}
\DeclareMathOperator{\Id}{Id}
\allowdisplaybreaks[2]
\title[The conformal Einstein-Dirac system on spheres]
{Coupling perturbations of the conformal Einstein-Dirac system on round spheres}
\author{Ke Xu }
\author{\quad Tian Xu}
\author{\quad Minbo Yang}
\address{Ke Xu \newline Department of Mathematics, \newline Zhejiang Normal University, \newline Jinhua, Zhejiang, 321004, China}
\email{xuke18264006093@163.com}
\address{Tian Xu \newline Department of Mathematics, \newline Zhejiang Normal University, \newline Jinhua, Zhejiang, 321004, China}
\email{xutian@amss.ac.cn}
\address{Minbo Yang \newline Department of Mathematics, \newline Zhejiang Normal University, \newline Jinhua, Zhejiang, 321004, China}
\email{mbyang@zjnu.edu.cn}
\date{}

\begin{document}
\begin{abstract}
	We study the conformal Einstein-Dirac system on the round
	sphere $S^m$, $m\geq3$, with coupling coefficient $1+\eps k$,
	where $k\in C^2(S^m,\R)$ and $\eps$ is small.
	Morse and degree conditions yield solutions with positive
	scalar component, and a generic nondegeneracy condition gives
	multiplicity modulo spinorial symmetries.
	For a suitable fixed coefficient, we obtain a noncompact
	solution sequence with infinitely many distinct energy values.
	A Kazdan-Warner identity gives nonexistence for a smooth
	Morse potential and its small $C^2$-perturbations, even with
	positive coupling.
	The existence proofs use a local reduction in the interior
	of the positive cone, accommodating the negative powers
	occurring when $m\geq5$.
\end{abstract}
\maketitle

\section{Introduction}

Let $(S^m,g)$ be the unit round sphere with its standard spin structure,
where $m\geq3$, and let $\Sigma S^m$ be the complex spinor bundle. Its
complex rank is $N=2^{\lfloor m/2\rfloor}$. We use the convention
$\Delta_g=\operatorname{div}_g\nabla$ and put
\[
 L_{g}=-\Delta_g+\frac{m-2}{4(m-1)}\Scal_g,
\]
for the Yamabe operator and write $D_{g}$ for its Dirac operator. More precisely, on $(S^m,g)$,
\begin{equation}\label{constants}
 L_{g}=-\Delta_g+\frac{m(m-2)}4.
\end{equation}
For $k\in C^2(S^m;\R)$ and a small real parameter $\eps$, we study
\begin{equation}\label{system}
 \left\{\begin{aligned}
 &L_{g}v=\dfrac{1+\eps k}{m-2}v^{(4-m)/(m-2)}|\psi|^2,\\[0.2em]
 &D_{g}\psi=(1+\eps k)v^{2/(m-2)}\psi,
 \end{aligned}\right.
 \qquad v>0,\quad\psi\in\Gamma(\Sigma S^m).
\end{equation}
Here $\Gamma(\Sigma S^m)$ denotes sections of the spinor bundle, with
regularity specified below. We use $\langle\ ,\ \rangle$ for the real
part of the Hermitian metric on $\Sigma S^m$ throughout. The Euler-Lagrange functional for \eqref{system} is
\begin{equation}\label{energy}
 E_\eps(v,\psi)=\frac12\int_{S^m}
 \big(vL_{g}v+\rp{D_{g}\psi}{\psi}-(1+\eps k)v^{2/(m-2)}|\psi|^2\big)\dd.
\end{equation}

For constant coupling, Maalaoui and Martino~\cite{MM19}
studied the conformal Dirac-Einstein system on closed
spin $3$-manifolds. They characterized Palais-Smale
sequences and obtained ground state solutions under an
Aubin-type strict inequality. This existence criterion
does not apply to the round sphere. In dimension four, Maalaoui and
Martino~\cite{MM26} established ground state existence for
the constant-coupling system on closed spin manifolds
whose conformal class has positive Yamabe constant.

For perturbations of the coupling coefficient on the round
$3$-sphere, existence results were obtained by Guidi,
Maalaoui, and Martino~\cite{GMM}. In the present paper,
we study \eqref{system} on the round sphere in every
dimension $m\geq3$, with particular emphasis on how the
potential $k$ determines existence, multiplicity, and
noncompactness of the solution set, as well as obstructions
to existence.

Our approach follows the finite-dimensional perturbation
scheme developed by Isobe and Xu~\cite{IX} for the
spinorial Yamabe equation, see also
the earlier work of Isobe~\cite{Isobe13}. Two features of the coupled problem
require a separate analysis. First, the kernel of the
linearized system contains both conformal and spinorial
directions, and its full real dimension must be computed.
Second, when $m\geq5$, the scalar equation contains a
negative power of $v$. The global $C^2$ Hilbert-space
formulation in \cite[Section~2]{IX} therefore cannot be
invoked directly.

The geometric meaning of the coefficient can be made explicit. For a
positive scalar component $v$, one can set
\[
 \widetilde g=v^{4/(m-2)}g,\qquad
 \widetilde\psi=v^{-(m-1)/(m-2)}\beta\psi,
\]
where $\beta:\Sigma_gS^m\to\Sigma_{\widetilde g}S^m$ is the canonical
fiberwise isometry. Conformal covariance transforms \eqref{system} into
\begin{equation}\label{geometric-system}
\left\{ \begin{aligned}
 &D_{\widetilde g}\widetilde\psi=(1+\eps k)\widetilde\psi,\\[0.2em]
 &\displaystyle\Scal_{\widetilde g}
 =\frac{4(m-1)}{(m-2)^2}(1+\eps k)|\widetilde\psi|_{\widetilde g}^2.
 \end{aligned}\right.
\end{equation}
Indeed, $L_{g}v=\frac{m-2}{4(m-1)}\Scal_{\widetilde g}
 v^{(m+2)/(m-2)}$ and
$|\psi|_g^2=v^{2(m-1)/(m-2)}|\widetilde\psi|_{\widetilde g}^2$;
substitution gives the second equation, and the Dirac transformation
law gives the first. Thus $1+\eps k$ prescribes the coefficient relating
scalar curvature to spinor density, as well as the zero-order term in
the transformed Dirac equation. In particular, our solutions determine metrics of positive scalar
curvature for sufficiently small $|\eps|$.

There is also a variational interpretation in the conformal
Einstein-Dirac setting. Under the same change of variables,
\begin{equation}\label{geometric-functional}
 E_\eps(v,\psi)=\frac12\int_{S^m}\bigg(
 \frac{m-2}{4(m-1)}\Scal_{\widetilde g}
 +\langle D_{\widetilde g}\widetilde\psi,\widetilde\psi\rangle
 -(1+\eps k)|\widetilde\psi|_{\widetilde g}^2
 \bigg)\,d\mathrm{vol}_{\widetilde g}.
\end{equation}
This combines total scalar curvature with a Dirac action and a
prescribed, spatially varying mass-type term. The connection between
gravity and spinor fields through the Einstein-Dirac functional has been
studied in \cite{KF}. Formula \eqref{geometric-functional}
provides the corresponding interpretation for our normalization in
arbitrary dimension. The metric is varied within its conformal class;
we do not assert the full tensor Einstein equation obtained from
unrestricted metric variations.

The difficulty caused by negative powers is local in the present
perturbation problem.
We work on $(S^m,g)$ in
\begin{equation}\label{positive-space}
 X=W^{2,p}(S^m;\R)\times W^{1,p}(\Sigma S^m),\qquad
 X_+=\big\{(v,\psi)\in X:\min_{S^m}v>0\big\},\quad p>m.
\end{equation}
The scalar nonnegative cone is
$\{v\in W^{2,p}(S^m;\R):v\geq0\text{ a.e.}\}$.
We identify each Sobolev class with its continuous representative,
which exists by the Sobolev embedding for $p>m$. Nonnegativity almost
everywhere is then equivalent to nonnegativity at every point.
The interior of this cone in the $W^{2,p}$ norm consists exactly
of the functions with $\min_{S^m} v>0$.
Accordingly, $X_+$ is the product of this interior with the spinor
space. The use of positive cones and ordered structures in the variational
study of positive and sign-changing elliptic solutions has a
substantial literature; see, for example,
\cite{Bartsch01,BLW, Hofer82, LiuSun01,LiWang02, ZhangLi03, ZouLi10} and references therein.
Here the interior of the scalar positive cone provides an open
domain for the local reduction. The $X$-norm is stronger than the natural
energy norm on $H^1(S^m)\times H^{1/2}(\Sigma S^m)$ and controls
uniform convergence. Here it serves a local
implicit-function argument.

Let $\Kcal$ be the set of real functions $k\in C^2(S^m,\R)$ satisfying
\begin{enumerate}
\item[(K1)] $k$ is Morse and $\Delta_g k(\xi)\neq0$ for every $\xi\in\Crit(k)$;
\item[(K2)] the following inequality holds:
\begin{equation}\label{degree-condition}
 S(k):=\sum_{\substack{\xi\in\Crit(k)\\ \Delta_g k(\xi)<0}}
 (-1)^{\ind(k,\xi)}\neq(-1)^m.
\end{equation}
\end{enumerate}
Here $\ind(k,\xi)$ is the Morse index of $k$ at $\xi$, namely the
number of negative eigenvalues of its Hessian.
Write $\omega_m=\Vol(S^m,g)$ and
\begin{equation}\label{cm}
 C_m=\frac{m^{m-1}(m-2)^2}{2^m}.
\end{equation}

\begin{theorem}[Existence and orbit counts]\label{existence}\label{multiplicity}
For every $m\geq3$ and $k\in\Kcal$, there exists $\eps_0>0$ such that
\eqref{system} has a classical solution $(v_\eps,\psi_\eps)$ for every
$0<|\eps|<\eps_0$. These solutions satisfy
\begin{equation}\label{positive-conclusion}
 v_\eps(x)>0,\qquad |\psi_\eps(x)|>0, \qquad x\in S^m
\end{equation}
and
\begin{equation}\label{energy-asymptotic}
 \int_{S^m}(1+\eps k)v_\eps^{2/(m-2)}|\psi_\eps|^2\dd
 =C_m\omega_m+O(|\eps|).
\end{equation}
Moreover, there is a residual subset $\Acal\subset C^2(S^m;\R)$ such that,
for $k\in\Acal\cap\Kcal$ and every sufficiently small nonzero $\eps$,
\eqref{system} has at least $N$ distinct $S^1$ solution orbits.
If $m\equiv2,3,4\pmod8$, it also has at least $N/2$ distinct
$S^3$ solution orbits. The solutions obtained satisfy
\eqref{positive-conclusion} and \eqref{energy-asymptotic}.
\end{theorem}

The equivalence relation in Theorem~\ref{multiplicity} is specified
by a fixed group action. The circle acts by
$(v,\psi)\mapsto(v,e^{i\theta}\psi)$, with $\theta\in\R$ constant
on $S^m$. In dimensions $m\equiv2,3,4\pmod8$, the action of unit quaternions commute with Clifford multiplication \cite[Section 1.7]{Friedrich} or \cite[Page 33, Table III]{LM}. Both actions preserve the equations and the functional. Two solutions belong to distinct $\mathsf G$-orbits if
there is no element $\tau\in\mathsf G$ taking one to the other,
where $\mathsf G=S^1$ or $S^3$, respectively. This is the precise
notion of distinctness used here. Clearly, orbit counts depend essentially on the symmetry group used to identify solutions. It does not require different
energies or different scalar components. 

Theorem \ref{multiplicity} doesn't assert inequivalence under every possible geometric symmetry. 
When one symmetry group is contained in another,
passing to the larger group may merge distinct orbits and
cannot increase their number. Distinguishing solutions up to
all natural symmetries of the fixed equation would therefore
require identifying its full symmetry group. In particular,
when $m=3$, we have $N=2$, and the lower bound for the number
of $S^3$-orbits in Theorem~\ref{multiplicity} is $N/2=1$.
Thus this theorem does not guarantee more than one
$S^3$-orbit in this dimension.

Additional geometric symmetries may arise from
orientation-preserving M\"obius transformations of $S^m$, together with their induced conformal
action on the scalar and spinor fields. Our purpose here is
to relate existence and multiplicity to information about
the potential $k$. We state the orbit counts for the specified
spinorial actions and leave a systematic study of the full
symmetry group and its effect on these counts for future work. The injectivity
argument establishing these orbit counts is given in Section~\ref{sec:perturbation}.

The next result indicates that the lower bound for the number
of solutions can be multiplied by separating the potentials in
stereographic coordinates. Let $\pi_{p_0}:S^m\setminus\{p_0\}\to\R^m$
be a stereographic projection.

\begin{theorem}[Separated potentials]\label{separated}
Let $l\geq2$ and $k_1,\ldots,k_l\in\Acal\cap\Kcal$ have a common
global minimum point $p_0$. Given $z_1,\ldots,z_l\in\R^m$, set
\begin{equation}\label{sum-potential}
 k(x)=
 \begin{cases}
 \displaystyle\sum_{i=1}^l k_i\!\left(\pi_{p_0}^{-1}
       (\pi_{p_0}(x)-z_i)\right),&x\neq p_0,\\[4pt]
 \displaystyle\sum_{i=1}^l k_i(p_0),&x=p_0.
 \end{cases}
\end{equation}
If $\min_{i\neq j}|z_i-z_j|$ is sufficiently large, then, for all
sufficiently small nonzero $\eps$, \eqref{system} has at least $lN$
distinct $S^1$ solution orbits. In dimensions $m\equiv2,3,4\pmod8$,
it also has at least $lN/2$ distinct $S^3$ solution orbits.
The conclusions \eqref{positive-conclusion} and
\eqref{energy-asymptotic} hold for the constructed solutions.
If, in addition, $k_i=t_i h$ for a single $h\in\Acal\cap\Kcal$
with a global minimum at $p_0$ and pairwise distinct positive numbers
$t_1,\ldots,t_l$, then solutions belonging to the $l$ parameter regions
can be constructed with disjoint energy ranges. In particular, the
system has at least $l$ solutions with pairwise different energies.
\end{theorem}

Each translated stereographic map extends smoothly across $p_0$,
so the function in \eqref{sum-potential} belongs to $C^2(S^m,\R)$.
The separation is fixed before $\eps_0$ is chosen. Each solution
constructed in this theorem is close to one standard bubble; the
integer $l$ counts different parameter regions. Spatial separation
alone does not force energy separation. The additional amplitude
condition distinguishes the first-order energy contributions from the
different regions. After the translations have been fixed, these
energy differences dominate the quadratic remainders for sufficiently
small nonzero $\eps$. The proof also gives a more general sufficient
condition in terms of the critical values of the first-order
perturbation term, introduced in Section~\ref{sec:existence}.

\begin{theorem}[Fixed-coefficient noncompactness]\label{fixed-noncompactness}
For every $m\geq3$ and $p_0\in S^m$, there is a nonnegative
function $k\in C^2(S^m,\R)$ satisfying
$k(p_0)=0$, $\nabla_g k(p_0)=0$, and $\nabla_g^2k(p_0)=0$, and there is
$\eps_0>0$ with the following property. For each fixed $\eps\in(-\eps_0,\eps_0)\setminus\{0\}$, the single system \eqref{system} has a sequence
of classical solutions $(v_{\eps,i},\psi_{\eps,i})$ satisfying
\eqref{positive-conclusion}, with both suprema tending to infinity.
The fields tend locally uniformly to zero away from $p_0$, and
\begin{equation}\label{fixed-concentration-measure}
 (1+\eps k)v_{\eps,i}^{2/(m-2)}|\psi_{\eps,i}|^2\,d\mathrm{vol}_g
 \rightharpoonup C_m\omega_m\delta_{p_0}.
\end{equation}
The sequence is bounded and weakly converges to zero in
$H^1(S^m)\times H^{1/2}(\Sigma S^m)$, with no strongly convergent
subsequence. Its energy values tend to
$C_m\omega_m/(2(m-2))$ and contain infinitely many distinct values.
\end{theorem}

The potential in Theorem~\ref{fixed-noncompactness} is constructed
as an infinite sum of translated compactly supported potentials
and is flat to second order at $p_0$. It therefore lies outside
the class $\Kcal$. Both $k$ and $\eps$ remain fixed along the
solution sequence. Uniform estimates after conformal recentering
allow the perturbation interval to be chosen independently of
the solution index.

We also derive a Kazdan-Warner obstruction for the coupled
system, see Proposition \ref{coupled-kw}. For the smooth height potential
$k(x)=x_{m+1}$ with $x\in S^m\hookrightarrow\R^{m+1}$, which satisfies \textup{(K1)} and fails
\textup{(K2)}, system \eqref{system} has no classical solution
with positive scalar component for any $\eps\neq0$.
This nonexistence conclusion persists throughout a sufficiently
small $C^2$ neighborhood of the height potential, with
$1+\eps k$ remaining positive whenever $|\eps|<1/2$.
Consequently, for each fixed sufficiently small nonzero $\eps$,
the set of potentials admitting such a solution is not dense
and hence cannot be residual in the whole space $C^2(S^m,\R)$.

The paper is organized as follows.
Section~\ref{sec:perturbation} presents the local reduction
and the finite-dimensional persistence arguments.
Section~\ref{sec:nondegeneracy} computes the full real kernel
of the coupled linearized system and verifies the analytic
hypotheses for reduction in $X_+$.
Section~\ref{sec:existence} analyzes the first-order perturbation
term and proves Theorems~\ref{existence} and~\ref{separated}.
Section~\ref{noncompactness} establishes fixed-coefficient
noncompactness through the translated-potential construction
and a uniform reduction after conformal recentering.
Finally, Section~\ref{sec:obstruction} derives the
Kazdan-Warner identity and proves nonexistence for the system \eqref{system}.

\section{A local perturbation framework}\label{sec:perturbation}

We adapt the variational perturbation framework of
\cite[Chapter~2]{AMbook} to functionals defined on an open
subset of a Banach space. In the application, this subset
is $X_+$, equipped with the $W^{2,p}\times W^{1,p}$ topology,
which is stronger than the natural $H^1\times H^{1/2}$ energy
topology. The finite-dimensional variational argument of
\cite[Section~2]{IX} is retained.

Let $X\hookrightarrow Y\hookrightarrow H$ be continuously embedded
real Banach spaces, where $H$ is a Hilbert space and $X$ is dense in $H$.
In this notation, $X$ is the space of fields and variations, $Y$ is the
target space for the equation, and $H$ supplies an inner product for
the first variation and the finite-dimensional orthogonal projections.
Let $U\subset X$ be open, let $E_0,G:U\to\R$ be $C^4$ functionals, and set
$E_\eps=E_0-\eps G$. The minus sign is only a convention: the case
$E_0+\eps G$ can be obtained by replacing $G$ with $-G$.

By definition, the Fr\'echet derivative of $E_\eps$, denoted by $dE_\eps(u)\in X^*$ (also denoted by $E_\varepsilon'(u)$ for notational convenience) is a bounded linear functional on $X$. An element $y\in Y$ may be used to define a functional by the rule
$h\mapsto(y,h)_H$, since
\[
 |(y,h)_H|\leq\|y\|_H\|h\|_H
 \leq C\|y\|_Y\|h\|_X\qquad \forall h\in X.
\]
The same rule for $y\in H$ gives a continuous injection of $H$
into $X^*$: if all these pairings vanish, density of $X$ in $H$
implies $y=0$. Thus the inner product gives the continuous inclusions
\[
 X\hookrightarrow Y\hookrightarrow H\hookrightarrow X^*.
\]
The last inclusion need not be surjective. In particular, neither
density nor differentiability of $E_\eps$ ensures that its derivative
has a representative in $Y$.

We impose this additional representation hypothesis: there is an
equation map $F_\eps:U\to Y$ such that
\begin{equation}\label{gradient-pairing}
 dE_\eps(u)[h]=(F_\eps(u),h)_H\qquad \forall h\in X.
\end{equation}
Here the right-hand side identifies $F_\eps(u)$ with an element of
$X^*$ through the preceding inclusion. Density guarantees uniqueness
of this representative, and hence
$dE_\eps(u)=0$ if and only if $F_\eps(u)=0$. This is not an application
of the Riesz representation theorem to an arbitrary element of $X^*$. The stronger assertion that the representative belongs to $Y$
records the regularity of the equation. In the Lyapunov-Schmidt reduction below, we shall formulate the normal equation and its linearization between the corresponding subspaces of $X$ and $Y$. We assume that $(\varepsilon,u)\mapsto F_\varepsilon(u)$ is jointly $C^3$ on $(-\varepsilon_1,\varepsilon_1)\times U$ for some $\varepsilon_1>0$.

Assume that $\Zcal\subset U$ is a finite-dimensional $C^4$ embedded
manifold with $F_0|_{\Zcal}=0$ and $E_0|_{\Zcal}\equiv\mathrm{const}$.
For $z\in\Zcal$, put $A_z=F_0'(z):X\to Y$ and assume
\begin{equation}\label{abstract-hypotheses}
 \ker A_z=T_z\Zcal,\qquad A_z\text{ is Fredholm of index zero},
 \qquad (A_zh_1,h_2)_H=(h_1,A_zh_2)_H.
\end{equation}
Let $P_z$ be the $H$-orthogonal projection onto $T_z\Zcal$.
It restricts to bounded projections on $X$ and $Y$. Indeed, choose a local
$C^3$ tangent frame $\vartheta_1(z),\ldots,\vartheta_d(z)$ and set
$M(z)_{ij}=(\vartheta_i(z),\vartheta_j(z))_H$, then
\[
 P_zh=\sum_{i,j=1}^d\vartheta_i(z)
       [M(z)^{-1}]_{ij}(h,\vartheta_j(z))_H \qquad \text{for } h\in X \text{ or } Y.
\]
The matrix is invertible because the frame is linearly independent.
Its entries and inverse are $C^3$ in $z$. Each pairing is continuous
on $X$ and $Y$, and the range lies in $T_z\Zcal\subset X$. Thus $P_z$ is
bounded on both spaces and depends $C^3$ on $z$.

\begin{proposition}[Local normal reduction]\label{normal-reduction}
On a neighborhood of each compact subset $\Zcal_c\subset\Zcal$,
there is a unique small normal correction
$w_\eps(z)\in(I-P_z)X$ satisfying
\begin{equation}\label{normal-equation}
 (I-P_z)F_\eps(z+w_\eps(z))=0,\qquad w_0(z)=0.
\end{equation}
It satisfies $\|w_\eps\|_{C^2(\Zcal_c;X)}=O(|\eps|)$.
The reduced functional
$E_\eps^{\mathrm{red}}(z)=E_\eps(z+w_\eps(z))$ has the expansion
\begin{equation}\label{reduced-expansion}
 E_\eps^{\mathrm{red}}
 =E_0|_{\Zcal}-\eps G|_{\Zcal}+\eps^2R_\eps,
\end{equation}
where $\{R_\eps:|\eps|<\eps_0\}$ is a bounded family in
$C^2(\Zcal_c,\R)$. 
Every critical point of $E_\eps^{\mathrm{red}}$ gives a critical
point of $E_\eps$.
\end{proposition}

\begin{proof}
Condition \eqref{abstract-hypotheses} gives $\Ran A_z\subset(I-P_z)Y$. These two closed subspaces have codimension $\dim T_z\Zcal$ by the index-zero hypothesis, so they
coincide in the sense that
\[
 A_z:(I-P_z)X\longrightarrow(I-P_z)Y
\]
is an isomorphism. In local trivializations the inverses vary
continuously with $z$ and are uniformly bounded on compact sets.
The parameter implicit function theorem gives a $C^3$ solution of
\eqref{normal-equation} near each $z$. 
Since $w_0(z)=0$, the mixed derivative
bounds through order three give the stated $O(|\eps|)$ estimate
in $C^2$. 

By \eqref{normal-equation}, $F_\eps(z+w_\eps(z))\in T_z\Zcal$.
For $\vartheta\in T_z\Zcal$,
\[
 dE_\eps^{\mathrm{red}}(z)[\vartheta]
 =(F_\eps(z+w_\eps),\vartheta+dw_\eps(z)[\vartheta])_H.
\]
The map $\vartheta\mapsto
P_z(\vartheta+dw_\eps(z)[\vartheta])$ on $T_z\Zcal$ is
$I+O(|\eps|)$ and hence invertible. It follows that
$dE_\eps^{\mathrm{red}}(z)=0$ implies $F_\eps(z+w_\eps(z))=0$.

To verify the remainder estimate, using Taylor's formula we have that
\[
 E_\eps^{\mathrm{red}}(z)-E_0(z)+\eps G(z)
 =\int_0^1(1-s)E_0''(z+sw_\eps)[w_\eps,w_\eps]\,ds-\eps\int_0^1G'(z+sw_\eps)[w_\eps]\,ds.
\]
Denoted by
\[
h_\eps(z)=\int_0^1\frac{\partial w_t(z)}{\partial t}\Big|_{t=s\eps}\,ds,
\]
use $w_0(z)=0$ and the fundamental
theorem of calculus, it holds that
\[
 h_\eps(z)=\left\{
 \aligned
 &\partial_\eps w_0(z), & & \eps=0 \\[0.3em]
 &\frac{w_\eps(z)}{\eps}, & & \eps\neq0.
 \endaligned
 \right.
\]
Moreover, the above mapping $(\eps,z)\mapsto h_\eps(z)$ is of class $C^2$: differentiating the integral
up to twice gives
\[
\partial_\eps h_\eps(z)=\int_0^1s\,\frac{\partial^2 w_t(z)}{\partial t^2}\Big|_{t=s\eps}\,ds,
\]
\[
\partial_\eps^2 h_\eps(z)=\int_0^1s^2\,\frac{\partial^3 w_t(z)}{\partial t^3}\Big|_{t=s\eps}\,ds,
\]
\[
\partial_z^2h_\eps(z)=\int_0^1\frac{\partial}{\partial t}\partial_z^2w_t(z)\Big|_{t=s\eps}\,ds
\]
and the mixed derivative can be treated similarly. These formulas involve at most the derivative of $w_\eps(z)$ up to order three, and this concludes the regularity.

Divide the Taylor identity above by $\eps^2$, we obtain
\[
R_\eps(z)=\int_0^1(1-s)E_0''(z+sw_\eps)[h_\eps,h_\eps]\,ds-\int_0^1G'(z+sw_\eps)[h_\eps]\,ds.
\]
The factors $E_0''$ and $G'$ are respectively $C^2$ and $C^3$ with $h_\eps$ being $C^2$, so the resulting expression for
$R_\eps$ is $C^2$, including at $\eps=0$. Therefore, $\{R_\eps\}$ is a bounded family in $C^2(\Zcal_c,\R)$.
\end{proof}

The regularity assumptions in Proposition~\ref{normal-reduction}
are sufficient for its $C^2$ remainder estimate, rather than a claim
of optimal regularity. A $C^1$ equation map already suffices for
the implicit-function construction of the normal correction;
in a Hilbert gradient setting this corresponds to a $C^2$ functional.
The expansion \eqref{reduced-expansion} additionally controls two
parameter derivatives of the remainder after division by $\eps^2$.
In the proof this requires two derivatives of $E_0''$, while the
$C^2$ extension of $w_\eps/\eps$ uses two derivatives of
$\partial_\eps w$. These requirements explain the $C^4$ functional
and $C^3$ equation-map hypotheses, respectively. Merely assuming
$E_\eps\in C^2(U,\R)$ does not by itself justify this conclusion.
With $C^3$ functionals and $C^2$ equation maps, the same integral
formula directly gives only a $C^1$ extension of $w_\eps/\eps$;
$C^{3,\alpha}$ functionals and $C^{2,\alpha}$ equation maps,
$0<\alpha<1$, similarly give only $C^{1,\alpha}$ by this formula.
Weaker remainder estimates or additional variational cancellations
may permit weaker assumptions, but would require a different proof.
The stated finite regularity is sufficient for all the estimates
used here; in our application the equation maps are in fact
$C^\infty$ in the fields on $X_+$, even with the coefficient $k\in C^2(S^m,\R)$.

For orbit counting, suppose that
$\Zcal\simeq\Pcal\times\mathcal S$, with $\Pcal\subset\R^d$ open
and $\mathcal S$ a closed manifold (compact without boundary). We write the
parameters as $(\varrho,\ga)$. In the application,
$\varrho=(\lambda,\xi)\in(0,\infty)\times\R^m$ is the dilation
and translation parameter and $\ga\in S^{2N-1}\subset\C^N$ is a
unit spinor parameter. Let a compact Lie group $\mathsf G$ act
continuously and linearly on $X,Y,H$, isometrically on $H$,
preserving $U$, the functional, and smoothly on
$\mathcal S$, with
\[
 \tau\cdot(\varrho,\ga)=(\varrho,\tau\cdot\ga).
\]
We assume that the action on $\mathcal S$ is free, meaning that
$\tau\cdot\ga=\ga$ implies $\tau=1$. The quotient
$Q=\mathcal S/\mathsf G$ then parametrizes the orbits and is a
compact manifold. Points $(\varrho,[\ga])$ of $\Pcal\times Q$
represent standard solutions up to exactly this group action.

Invariance of the functional, uniqueness of the representative in
\eqref{gradient-pairing}, and the isometric action on $H$ imply
equivariance of the equation. The orthogonal projections are equivariant
because the action preserves the tangent spaces.
Uniqueness in \eqref{normal-equation} gives
\begin{equation}\label{equivariance}
 w_\eps(\tau\cdot z)=\tau\cdot w_\eps(z).
\end{equation}
Thus $E_\eps^{\mathrm{red}}$ is constant on group orbits and defines
a well-defined function on $\Pcal\times Q$.

\begin{corollary}[Local multiplicity]
\label{abstract-multiple}
Suppose that the restriction $G|_{\Zcal}=G(\varrho,\ga)=\mathcal T(\varrho)$
is independent of $\ga$, and that $\varrho_0$ is a nondegenerate
critical point of $\mathcal T$. For every sufficiently small nonzero
$\eps$, there are at least $\Cat(Q)$ distinct $\mathsf G$ orbits of
critical points of $E_\eps$ near $\{\varrho_0\}\times\mathcal S$.
\end{corollary}

Here $\Cat(Q)$ denotes the Lusternik-Schnirelmann category
of $Q$, namely the least number of open subsets of $Q$,
each contractible in $Q$, required to cover $Q$.

\begin{proof}
By Proposition~\ref{normal-reduction},
$(E_0|_{\Zcal}-E_\eps^{\mathrm{red}})/\eps$ extends to a $C^2$
function of $(\eps,\varrho,\ga)$ and equals $\mathcal T$ at
$\eps=0$. Apply the implicit function theorem to its derivative
in $\varrho$ at $\varrho_0$. Compactness of $\mathcal S$ and
uniqueness provide a single $C^1$ branch
\[
 \varrho_\eps(\ga)=\varrho_0+O(|\eps|)
 \]
 such that
 \[
\varrho_\eps(\tau\cdot\ga)=\varrho_\eps(\ga) \quad \text{and} \quad \partial_\varrho E_\eps^{\mathrm{red}}(\varrho_\eps(\ga),\ga)=0.
\]
Hence
$\ga\mapsto E_\eps^{\mathrm{red}}(\varrho_\eps(\ga),\ga)$ defines
a $C^1$ function on $Q$. The Lusternik--Schnirelmann theorem
\cite[Chapter~II, Section~3.2]{Chang93} gives at least
$\Cat(Q)$ distinct critical points on $Q$. 

It remains to check that distinct quotient parameters cannot
produce the same solution orbit. Consider the normal-coordinate map
\[
 (z,w)\longmapsto z+w,\qquad w\in(I-P_z)X.
\]
Its differential at $(z,0)$ is the isomorphism
$(\vartheta,\eta)\mapsto\vartheta+\eta$ from
$T_z\Zcal\oplus(I-P_z)X$ onto $X$. The inverse function theorem
therefore gives unique normal coordinates near each $(z,0)$.
This uniqueness is uniform for sufficiently small normal vectors
with base in a fixed compact subset: otherwise one could find
$z_j+w_j=z'_j+w'_j$ with distinct normal coordinates, compact base
points, and $\|w_j\|_X+\|w'_j\|_X\to0$. Passing to a subsequence,
$z_j\to z$ and $z'_j\to z'$. Equality of the images implies
$z=z'$, so both pairs eventually lie in the same local coordinate
neighborhood, a contradiction.

Now write $u_\eps(z)=z+w_\eps(z)$. If
$u_\eps(z')=\tau\cdot u_\eps(z)$, equivariance gives
$u_\eps(z')=u_\eps(\tau\cdot z)$. Apply the preceding uniqueness
on a compact invariant parameter neighborhood. It follows that
$z'=\tau\cdot z$. Thus equal solution orbits would imply equal
quotient parameters. This proves the claimed distinction of orbits.
\end{proof}

Nonzero degree can give existence without an individually
nondegenerate critical point. We use the following additional argument.

\begin{lemma}[Persistence on a compact quotient]\label{index-lemma}
Let $\Omega\Subset\R^d$ have smooth boundary, let $Q$ be a closed
oriented manifold with Euler characteristic $\chi(Q)\neq0$, and let $\mathcal T$ be smooth
near $\overline\Omega$. Suppose that
\[
 \nabla\mathcal T\neq0\text{ on }\partial\Omega,\qquad
 \deg(\nabla\mathcal T,\Omega,0)\neq0.
\]
If $f_\eps(\varrho,\zeta)=\mathcal T(\varrho)+O(|\eps|)$ in
$C^1(\overline\Omega\times Q)$, where $\zeta\in Q$, then
$f_\eps$ has an interior critical point for small $|\eps|$.
\end{lemma}

\begin{proof}
The vector fields $\nabla f_\eps$ and
$(\nabla\mathcal T,0)$ are homotopic without zeros on
$\partial\Omega\times Q$, since their base components are close
in the uniform norm. For a field nonvanishing on the boundary,
the sum of the local indices of a transverse perturbation is
independent of that perturbation and invariant under such a
homotopy. 

Choose a sufficiently small regular value $b$ of $\nabla\mathcal T$
and a Morse function $h$ on $Q$. The field
\[
 (\nabla\mathcal T(\varrho)-b,\delta\nabla h(\zeta)),\qquad\delta>0,
\]
is homotopic to $(\nabla\mathcal T,0)$ without boundary zeros.
Its derivative at each zero is block diagonal, so its sum of
local indices equals
\[
 \deg(\nabla\mathcal T,\Omega,0)
 \sum_{\zeta\in\Crit(h)}(-1)^{\ind(h,\zeta)}
 =\deg(\nabla\mathcal T,\Omega,0)\chi(Q)\neq0.
\]
The equality involving $\chi(Q)$ follows from the Morse
relation for $h$; see \cite[Chapter~I, Section~4]{Chang93}. And since vector field without zeros has index sum zero, we find $\nabla f_\eps$ must vanish in the interior.
\end{proof}

\begin{remark}
Corollary~\ref{abstract-multiple} is the local version of
\cite[Theorem 2.4 and Remark 2.5]{IX} after taking the specified
symmetry quotient. Lemma~\ref{index-lemma} proves existence even
when the base critical points are degenerate. It alone gives no
multiplicity bound: a degenerate isolated zero need not contribute
$+1$ or $-1$, so a nonzero sum of indices cannot in general be
interpreted as that many distinct zeros. The orbit count instead
uses category and the nondegenerate base point.
\end{remark}

\section{The standard family and its normal nondegeneracy}\label{sec:nondegeneracy}

\subsection{Conformal parametrization}
Ground state solutions with positive scalar component for
the three-dimensional conformal Dirac-Einstein system
were classified by Borrelli and Maalaoui~\cite{BMa21}.
Their scalar and spinorial components are given by
Aubin-Talenti profiles and conformal images of Killing
spinors, respectively. For the critical Dirac equation,
the corresponding classification is established in
\cite{BMW21}. We now construct a standard solution family
for \eqref{system} with $\eps=0$ in every dimension
$m\geq3$ and determine its normalization. Let
\begin{equation}\label{ab}
 A=\left(\frac m2\right)^{(m-2)/2},\qquad B=\frac{m-2}{\sqrt2}A.
\end{equation}
Choose a unit Killing spinor $\Xi$ satisfying
$\nabla_X\Xi=-\tfrac12X\cdot\Xi$. Then $D_{g}\Xi=\frac m2\Xi$ and
$|\Xi|=1$. The identities
\begin{equation}\label{normalization}
 A^{2/(m-2)}=\frac m2,\qquad \frac1{m-2}A^{(6-2m)/(m-2)}B^2=\frac{m(m-2)}4,
 \qquad A^{2/(m-2)}B^2=C_m
\end{equation}
show that
$(A,B\Xi)\in\big(C^\infty(S^m;\R)\times C^\infty(\Sigma S^m)\big)
\cap X_+$ solves \eqref{system} with $\eps=0$.

In stereographic coordinates, $g=f_0^2g_{\R^m}$ with
$f_0(x)=2/(1+|x|^2)$. Write $\beta$ for the conformal identification
of the spinor bundles. If
\[
 U=f_0^{(m-2)/2}(v\circ\pi^{-1}),\qquad
 \Phi=f_0^{(m-1)/2}\beta^{-1}(\psi\circ\pi^{-1}),
\]
then conformal covariance gives
\begin{align*}
 (L_{g}v)\circ\pi^{-1}&=f_0^{-(m+2)/2}(L_{g_{\R^m}}U),\\
 (D_{g}\psi)\circ\pi^{-1}&=f_0^{-(m+1)/2}\beta(D_{g_{\R^m}}\Phi).
\end{align*}
The energy and the integral $\int v^{2/(m-2)}|\psi|^2\dd$ are preserved.
In particular, \eqref{system} becomes the corresponding Euclidean
system
\begin{equation}\label{euclidean-system}
\left\{\begin{aligned}
 &L_{g_{\R^m}}U=\dfrac{1+\eps k\circ\pi^{-1}}{m-2}
 U^{(4-m)/(m-2)}|\Phi|^2,\\[0.2em]
 &D_{g_{\R^m}}\Phi=(1+\eps k\circ\pi^{-1})U^{2/(m-2)}\Phi,
 \end{aligned}\right.
\end{equation}
where $L_{g_{\R^m}}=-\Delta_{\R^m}$ and $D_{g_{\R^m}}$ denote
only the Euclidean operators. The subscripts will be retained.

Set $\Pcal=(0,\infty)\times\R^m$. For $\varrho=(\lambda,\xi)\in\Pcal$
and a unit spinor $\ga\in\C^N$, define
\begin{equation}\label{bubbles}
 \begin{split}
 f_{\lambda,\xi}(x)&=\frac{2\lambda}{\lambda^2+|x-\xi|^2},
 \qquad U_{\lambda,\xi}=Af_{\lambda,\xi}^{(m-2)/2},\\
 \Phi_{\lambda,\xi,\ga}(x)&=
 B\,2^{(m-1)/2}\frac{\lambda^{(m-1)/2}}
 {(\lambda^2+|x-\xi|^2)^{m/2}}
 \big(\lambda-(x-\xi)\big)\cdot\ga.
 \end{split}
\end{equation}
Here Clifford multiplication is extended to scalar-plus-vector
expressions by $(\lambda-y)\cdot\ga=\lambda\ga-y\cdot\ga$.
For a tangent vector $V$, it satisfies
$V\cdot V\cdot\ga=-|V|^2\ga$ and is skew-adjoint. Hence
\begin{equation}\label{bubble-identities}
 |\Phi|=Bf_{\lambda,\xi}^{(m-1)/2},\qquad
 D_{g_{\R^m}}\Phi=\frac m2 f_{\lambda,\xi}\Phi,\qquad
 L_{g_{\R^m}}U=\frac{m(m-2)}4 A f_{\lambda,\xi}^{(m+2)/2}.
\end{equation}

Pulling the pairs $(U_{\lambda,\xi}, \Phi_{\lambda,\xi,\ga})$ back to $S^m$ gives a family
\begin{equation}\label{critical-family}
 \Zcal=\big\{z_{\lambda,\xi,\ga}=(v_{\lambda,\xi},\psi_{\lambda,\xi,\ga})\big\}\simeq\Pcal\times S^{2N-1},
 \qquad \dim_\R\Zcal=m+2N.
\end{equation}
All its scalar components and spinor norms are strictly positive.
The family is smooth on the whole sphere: it is obtained by
conformal transformations of the constant scalar and Killing spinor
solutions. Its scalar component determines $(\lambda,\xi)$, since
\[
 (U_{\lambda,\xi}/A)^{-2/(m-2)}
 =\frac{\lambda^2+|x-\xi|^2}{2\lambda}.
\]
For fixed $(\lambda,\xi)$ the spinor determines $\ga$. The differential
is injective by the same formulas, which also give local continuous
inverse parameters. Thus \eqref{critical-family} is an embedded
manifold in $X$.

The scaling formulas are
\begin{equation}\label{scaling}
 \begin{split}
 U_{\lambda,\xi}(x)&=\lambda^{-(m-2)/2}
 U_{1,0}\!\left(\frac{x-\xi}{\lambda}\right),\\
 \Phi_{\lambda,\xi,\ga}(x)&=\lambda^{-(m-1)/2}
 \Phi_{1,0,\ga}\!\left(\frac{x-\xi}{\lambda}\right).
 \end{split}
\end{equation}
In particular $|\Phi_{\lambda,\xi,\ga}|^2$ scales as $\lambda^{1-m}$.

\subsection{The full real kernel}

For $(v,\psi)\in X_+$, define the operator representing the first
variation of $E_\eps$ by
\begin{equation}\label{euler-operator}
 F_\eps(v,\psi)=
 \left(L_{g}v-\frac1{m-2}(1+\eps k)v^{(4-m)/(m-2)}|\psi|^2,
             D_{g}\psi-(1+\eps k)v^{2/(m-2)}\psi\right).
\end{equation}
Thus $F_0$ is obtained by setting $\eps=0$; its zeros solve the
unperturbed system. The mapping properties will be verified in Section~\ref{subsec:mapping property}.

\begin{proposition}\label{kernel}
Let $A_z=F_0'(z)$. Then
\[
 \ker A_z=T_z\Zcal,\qquad \dim_\R\ker A_z=m+2N
 \quad(z\in\Zcal).
\]
\end{proposition}

\begin{proof}
We first work at $(A,B\Xi)$. Under the change of unknowns $v=Au$ and $\psi=B\phi$,
the unperturbed system is equivalent to
\[
 L_{g}u=\frac{m(m-2)}4 u^{(4-m)/(m-2)}|\phi|^2,\qquad D_{g}\phi=\frac m2 u^{2/(m-2)}\phi.
\]
The standard solution $(A,B\Xi)$ corresponds to $(u,\phi)=(1,\Xi)$.
At $(1,\Xi)$ let $(s,\eta)$ be a real variation and put
$r=\rp{\eta}{\Xi}$. The linearized equations are
\begin{equation}\label{linear-system}
 \left(-\Delta_g+\frac{m(m-3)}2\right)s
 =\frac{m(m-2)}2r,\qquad
 \left(D_{g}-\frac m2\right)\eta=\frac m{m-2}s\Xi.
\end{equation}
For spinor fields we use the connection Laplacian, with the sign
convention
\[
 \Delta_g^\Sigma\eta=-\nabla^*\nabla\eta
 =\sum_{j=1}^m\left(\nabla_{e_j}\nabla_{e_j}\eta
                  -\nabla_{\nabla^g_{e_j}e_j}\eta\right).
\]
Here $\nabla$ is the spin connection and $\nabla^g$ the Levi-Civita
connection. This is the bundle analogue of the scalar
Laplace-Beltrami operator $\Delta_g$; it acts on spinors, not on
functions. With this convention the Schr\"odinger-Lichnerowicz
formula is $D_{g}^2=-\Delta_g^\Sigma+\Scal_g/4$; see \cite{Friedrich}.
The product rule and the Killing equation give
\begin{align*}
 \Delta_g r
 &=\rp{\Delta_g^\Sigma\eta}{\Xi}
   +\rp{\eta}{\Delta_g^\Sigma\Xi}
   +2\sum_j\rp{\nabla_{e_j}\eta}{\nabla_{e_j}\Xi}\\
 &=-\rp{D_{g}^2\eta}{\Xi}+\frac{m(m-2)}4r+\rp{D_{g}\eta}{\Xi}.
\end{align*}
Here we used $D_{g}^2=-\Delta_g^\Sigma+m(m-1)/4$ and
$\Delta_g^\Sigma\Xi=-m\Xi/4$. From \eqref{linear-system},
\[
 D_{g}^2\eta=\frac{m^2}4\eta+\frac{m^2}{m-2}s\Xi+\frac m{m-2}\nabla s\cdot\Xi,
 \qquad \rp{X\cdot\Xi}{\Xi}=0.
\]
It follows that
\begin{equation}\label{r-laplacian}
 -\Delta_g r=\frac{m(m-1)}{m-2}s.
\end{equation}
Eliminating $r$ yields
\begin{equation}\label{factorization}
 (-\Delta_g-m)\left(-\Delta_g+\frac{m(m-1)}2\right)s=0.
\end{equation}
The two operators in parentheses on the left-hand side are polynomials in $-\Delta_g$, and hence commute. Thus we conclude that
\[
 s\in E_1:=\ker(-\Delta_g-m),\qquad
 r=\frac{m-1}{m-2}s,\qquad\dim_\R E_1=m+1.
\]
For $s\in E_1$ define
\begin{equation}\label{particular-spinor}
 \eta_s=\frac{m-1}{m-2}s\Xi+\frac1{m-2}\nabla s\cdot\Xi.
\end{equation}
The identities
\[
 D_{g}(s\Xi)=\nabla s\cdot\Xi+\frac m2 s\Xi,\qquad
 D_{g}(\nabla s\cdot\Xi)=ms\Xi-\frac{m-2}{2}\nabla s\cdot\Xi
\]
show that $(s,\eta_s)$ satisfies \eqref{linear-system}. Let $(s,\eta)$ be any solution of \eqref{linear-system}.
The preceding argument shows that $s\in E_1$ and
$\rp{\eta}{\Xi}=\frac{m-1}{m-2}s$. Setting
$\zeta=\eta-\eta_s$, we obtain
\[
\left(D_g-\frac m2\right)\zeta=0,\qquad
\rp{\zeta}{\Xi}
=\rp{\eta}{\Xi}-\rp{\eta_s}{\Xi}=0.
\]
Consequently, the complete solution space of the linearized system \eqref{linear-system} is
\begin{equation}\label{full-kernel}
\Big\{(s,\eta_s+\zeta):s\in E_1,\ D_{g}\zeta=\frac m2\zeta,
                    \ \rp{\zeta}{\Xi}=0 \Big\}.
\end{equation}

Every eigenspinor with eigenvalue $m/2$ is a Killing spinor with Killing constant $-1/2$. Indeed,
the integrated Lichnerowicz identity gives
\[
 \int_{S^m}\sum_j\left|\nabla_{e_j}\zeta+
                 \tfrac12e_j\cdot\zeta\right|^2\dd
 =\|\nabla\zeta\|_2^2-\frac m2\|\zeta\|_2^2
                    +\frac m4\|\zeta\|_2^2=0.
\]
And to see the space of such Killing spinors has complex dimension $N$, we can use the observation that evaluation at a point is injective by the Killing equation, while
\eqref{bubbles} at $(\lambda,\xi)=(1,0)$ supplies $N$ independent
complex directions. For two Killing spinors $\zeta$ and $\Xi$ with same Killing constant $-1/2$,
the skew-adjointness of Clifford multiplication gives
\[
X\rp{\zeta}{\Xi}
=-\frac12\big(
\rp{X\cdot\zeta}{\Xi}
+\rp{\zeta}{X\cdot\Xi}
\big)=0.
\]
Thus $\rp{\zeta}{\Xi}$ is constant on $S^m$, and its vanishing
is equivalent to its vanishing at an arbitrary given point $p\in S^m$. The map
\[
\zeta\longmapsto\rp{\zeta(p)}{\Xi(p)}
\]
is a nonzero real linear functional on the space of such
Killing spinors, since it takes the value $1$ at $\zeta=\Xi$. Consequently, the subspace satisfying
$\rp{\zeta}{\Xi}=0$ has real dimension $2N-1$. And therefore \eqref{full-kernel}
has real dimension $(m+1)+(2N-1)=m+2N$.

Differentiating $F_0(z_{\lambda,\xi,\ga})=(0,0)$ along the parameters
gives $T_z\Zcal\subset\ker A_z$. We first prove equality at a
constant solution $z=z_{1,0,\ga}=(A,B\Xi)$. The preceding kernel
calculation uses the normalized variables; the invertible rescaling
$(s,\eta)\mapsto(As,B\eta)$ identifies its solution space with
$\ker A_z$.

The scalar component of the standard family satisfies
\[
v_{\lambda,\xi}\circ\pi^{-1}
=A\left(\frac{f_{\lambda,\xi}}{f_0}\right)^{(m-2)/2}.
\]
Since $f_{1,0}=f_0$, its derivatives at $(1,0)$ with respect to
$\lambda$ and $\xi_j$ are $A(m-2)/2$ times the corresponding
derivatives of $\log f_{\lambda,\xi}$. Direct calculation gives
\[
\partial_\lambda\log f_{\lambda,\xi}|_{(1,0)}
=\frac{|x|^2-1}{1+|x|^2},\qquad
\partial_{\xi_j}\log f_{\lambda,\xi}|_{(1,0)}
=\frac{2x_j}{1+|x|^2}.
\]
The functions on the right are the coordinate functions of $S^m$
in stereographic coordinates and form a basis of $E_1$. Thus the
$m+1$ conformal parameter directions are linearly independent,
and their scalar components span $E_1$.

At $(\lambda,\xi)=(1,0)$, the spinor formula is complex linear
in $\ga$ and, after division by $B$, identifies $\C^N$ with the
space of $-\frac12$-Killing spinors, preserving the pointwise
norm. Varying $\ga\in S^{2N-1}$ therefore gives precisely the
tangent vectors
\[
(0,B\zeta),\qquad D_g\zeta=\frac m2\zeta,\qquad
\rp{\zeta}{\Xi}=0.
\]
These directions have real dimension $2N-1$. Their scalar
components vanish, so they are independent of the conformal
directions. Consequently,
\[
\dim_\R T_z\Zcal=(m+1)+(2N-1)=m+2N=\dim_\R\ker A_z,
\]
and the inclusion $T_z\Zcal\subset\ker A_z$ is actually equality. 

Each $z_{\lambda,\xi,\ga}\in\Zcal$ is the conformal image of
$z_{1,0,\ga}$ under the transformation induced by Euclidean
translation and dilation. The induced maps on the field and
equation spaces are invertible. Differentiating the covariance
identity shows that the map on fields carries the linearized
kernel at $z_{1,0,\ga}$ onto that at $z_{\lambda,\xi,\ga}$.
Since this conformal transformation preserves $\Zcal$, the same
map carries the corresponding tangent spaces onto each other.
Hence $T_z\Zcal=\ker A_z$ for every $z\in\Zcal$.
\end{proof}


\subsection{Reduction in the interior of the positive cone}\label{subsec:mapping property}

Take $X,X_+$ from \eqref{positive-space} and let
\[
 Y=L^p(S^m;\R)\times L^p(\Sigma S^m),\qquad
 H=L^2(S^m;\R)\times L^2(\Sigma S^m),
\]
with the real $L^2$-inner product. These spaces have different roles: $X$
is the domain of the differential operators and controls the pointwise
positivity of the scalar component; $Y$ is the target space of the Euler-Lagrange operator $F_\eps$; $H$ provides the pairing representing the first variation and the finite-rank
orthogonal projections onto $T_z\Zcal$. In particular, $H$ here is an
auxiliary $L^2$ space, not the natural $H^1\times H^{1/2}$ energy space.
The mapping properties
\[
 L_g:W^{2,p}(S^m)\to L^p(S^m),\qquad
 D_g:W^{1,p}(\Sigma S^m)\to L^p(\Sigma S^m)
\]
explain this choice of $Y$. Elliptic estimates with right-hand sides
in $Y$ give an inverse with values in $X$. Merely taking $H$ as the
target would not give this regularity: $L^2$ right-hand sides need not
produce fields in $W^{2,p}\times W^{1,p}$.
Since $p>m\geq3$ and the sphere has finite
volume, $Y\hookrightarrow H$; smooth fields belong to $X$ and are dense
in $H$. For $u=(v,\psi)\in X_+$ and $h=(\sigma,\eta)\in X$,
integration by parts for $L_g$ and formal self-adjointness of $D_g$ give
\[
 dE_\eps(v,\psi)[\sigma,\eta]
 =\int_{S^m}\big[(F_\eps(v,\psi))_1\sigma
       +\rp{(F_\eps(v,\psi))_2}{\eta}\big]\dd,
\]
where $F_\eps(v,\psi)=((F_\eps(v,\psi))_1,(F_\eps(v,\psi))_2)$ is the operator in \eqref{euler-operator}.
Thus its components in $L^p$ characterize the representative of the derivative
in \eqref{gradient-pairing}.
Sobolev embedding gives $X\hookrightarrow C^0\times C^0$.
On bounded neighborhoods with $v\geq\delta>0$, every derivative
of a real power of $v$ has bounded coefficients. Multiplication
is controlled by $W^{1,p}\hookrightarrow L^\infty$.
It follows that $F_\eps:X_+\to Y$ and $E_\eps:X_+\to\R$ are
smooth. This smoothness concerns the field variables; $k\in C^2$
is more than sufficient for these maps.

All open sets and normal corrections in this section use the norm
$\|(v,\psi)\|_X=\|v\|_{W^{2,p}}+\|\psi\|_{W^{1,p}}$.
The embedding into $C^0\times C^0$ shows that a sufficiently small
$X$-ball about a point with $\min v>0$ remains in $X_+$.
Conversely, a nonnegative $v$ with a zero can be made negative there
by subtracting an arbitrarily small positive constant. This justifies
the asserted description of the interior of the scalar positive cone.
The situation is different in the natural energy topology. Around a smooth positive scalar
field, an $H^1$-ball contains functions that change sign.
Thus merely restricting the energy to the positive cone in
$H^1\times H^{1/2}$ does not provide an open neighborhood for
the implicit function theorem. Moreover, when $m\geq4$, $|v|^{2/(m-2)}$ is not differentiable at $v=0$. 

At $z=(v_z,\psi_z)\in\Zcal$, the derivative of the Euler-Lagrange operator gives
\begin{equation}\label{hessian-operator}
 \begin{split}
 (A_z(\sigma,\eta))_1
 &=L_{g}\sigma-\frac{4-m}{(m-2)^2}v_z^{(6-2m)/(m-2)}|\psi_z|^2\sigma
          -\frac2{m-2}v_z^{(4-m)/(m-2)}\rp{\psi_z}{\eta},\\
 (A_z(\sigma,\eta))_2
 &=D_{g}\eta-v_z^{2/(m-2)}\eta-\frac2{m-2}v_z^{(4-m)/(m-2)}\sigma\psi_z.
 \end{split}
\end{equation}
In particular $A_z$ is symmetric under the real $L^2$ pairing.
Equivalently, for $h_i=(\sigma_i,\eta_i)$,
\begin{align}\label{hessian-form}
 E_0''(z)[h_1,h_2]
 &=\int_{S^m}\Big[\sigma_1L_{g}\sigma_2+\rp{D_{g}\eta_1}{\eta_2}
       -v_z^{2/(m-2)}\rp{\eta_1}{\eta_2}\notag\\
 &\qquad -\frac{4-m}{(m-2)^2}v_z^{(6-2m)/(m-2)}|\psi_z|^2\sigma_1\sigma_2
       \notag\\
 &\qquad -\frac2{m-2}v_z^{(4-m)/(m-2)}\big(\sigma_1\rp{\psi_z}{\eta_2}
                    +\sigma_2\rp{\psi_z}{\eta_1}\big)\Big]\dd.
\end{align}

The diagonal map $\operatorname{diag}(L_{g},D_{g}):X\to Y$ is an
isomorphism by elliptic theory: $L_{g}>0$ and the Lichnerowicz formula
gives $\ker D_{g}=0$ on $(S^m,g)$. All remaining terms in
\eqref{hessian-operator} are compact from $X$ to $Y$, by the
compact embeddings into $L^p$. Hence $A_z$ is Fredholm of index
zero. Kernel elements are smooth by elliptic bootstrapping, so
Proposition~\ref{kernel} applies. All hypotheses of
Proposition~\ref{normal-reduction} are now verified.

For later use define
\begin{equation}\label{perturbation}
 G_k(v,\psi)=\frac12\int_{S^m}kv^{2/(m-2)}|\psi|^2\dd,
 \qquad E_\eps=E_0-\eps G_k.
\end{equation}
On every compact parameter domain the reduction gives
\begin{equation}\label{ed-reduction}
 \begin{split}
 E_\eps^{\mathrm{red}}(z)&=E_0(z)-\eps G_k(z)+O(\eps^2)
 \quad\text{in }C^2,\\
 E_0(z)&=\frac{C_m\omega_m}{2(m-2)}\qquad \forall z\in\Zcal.
 \end{split}
\end{equation}
The $S^1$ action $(v,\psi)\mapsto(v,e^{i\theta}\psi)$ and the $S^3$ action in dimensions $m\equiv2,3,4\pmod8$ preserve the equations,
the functional and $\Zcal$. Their action changes $\ga$ and fixes $(\lambda,\xi)$.
The actions on the unit spinor sphere are free. For the circle,
$e^{i\theta}\ga=\ga$ and $|\ga|=1$ imply $e^{i\theta}=1$.
For $S^3$, identify the spinor module with $\HH^{N/2}$; a nonzero
coordinate of $\ga$ shows that $\tau\ga=\ga$ implies $\tau=1$.
Thus the quotient construction of Section~\ref{sec:perturbation} applies to both groups.

\section{The first-order perturbation term and the main results}\label{sec:existence}

\subsection{Boundary behavior and topological degree}

For $k\in\Kcal$, choose a global minimum point $p_0$ as the
stereographic pole and put $a=k\circ\pi_{p_0}^{-1}$.
By \eqref{bubble-identities} and conformal invariance,
\begin{equation}\label{leading}
 G_k(z_{\lambda,\xi,\ga})=\frac{C_m}{2}T_k(\lambda,\xi),\qquad
 T_k(\lambda,\xi)=\int_{\R^m}a(\lambda y+\xi)\rho(y)\,dy,
 \quad \rho(y)=\frac{2^m}{(1+|y|^2)^m}.
\end{equation}
The right side is independent of the unit spinor parameter $\ga$.
The function $-(C_m/2)T_k$ is precisely the coefficient of $\eps$
in the expansion of $E_\eps^{\mathrm{red}}$ on $\Zcal$.
Thus $T_k$ determines the first-order finite-dimensional problem. Related applications of finite-dimensional reduction and
Morse theory to the prescribed scalar curvature problem
can be found in \cite{AGP99,Mal02}. The degree calculation
below uses an even extension in the scale parameter, as in
\cite{AGP99} and \cite[Section~3]{IX}. We give the details
for the function $T_k$ defined above.

\begin{lemma}\label{degree-lemma}
There is a smooth bounded open set $\Omega\Subset\Pcal$ containing
all critical points of $T_k$ in $\Pcal$, such that
\begin{equation}\label{degree-formula}
 \nabla T_k\neq0\text{ on }\partial\Omega,\qquad
 \deg(\nabla T_k,\Omega,0)=S(k)-(-1)^m\neq0,
\end{equation}
where $S(k)$ is given in \eqref{degree-condition}.
\end{lemma}

\begin{proof}
The function $a$ and its derivatives up to second order are bounded. Since $m\geq3$,
\begin{equation}\label{moments}
 \int_{\R^m}\rho=\omega_m,\qquad
 \int_{\R^m}y_iy_j\rho=\frac{\omega_m}{m-2}\delta_{ij}.
\end{equation}
These formulas follow from radial integration (or the beta
integral). It follows from the expression of $T_k$ in \eqref{leading} that the integrand and its parameter derivatives up to second order are
dominated by the integrable function $(1+|y|^2)\rho(y)$. Thus $T_k$
has an even $C^2$ extension $\widetilde T_k$ to
$\R\times\R^m$, and
\begin{equation}\label{boundary-hessian}
 \widetilde T_k(0,\xi)=\omega_m a(\xi),\qquad
 \nabla^2\widetilde T_k(0,\xi)=
 \begin{pmatrix}
 \dfrac{\omega_m}{m-2}\Delta_{\R^m} a(\xi)&0\\
 0&\omega_m \nabla^2a(\xi)
 \end{pmatrix}.
\end{equation}
In particular the critical points on $\lambda=0$ are precisely
$(0,\xi)$ with $\xi\in\Crit(a)$. They are nondegenerate and have
local index
\begin{equation}\label{boundary-index}
 \ind(\nabla\widetilde T_k,(0,\xi))
 =\operatorname{sgn}(\Delta_{\R^m} a(\xi))(-1)^{\ind(a,\xi)}.
\end{equation}
At a critical point, the spherical and Euclidean Hessians have
the same property, and $\Delta_{\R^m} a=f_0^2(\Delta_gk)\circ\pi^{-1}$.

We next analyze the behavior of $\nabla\widetilde T_k$ as
$\lambda^2+|\xi|^2\to\infty$.
Write $b(y)=\nabla a(y)\cdot y$. For $w\neq0$ near the origin, define
$\widehat a(w)=a(w/|w|^2)$ and set $\widehat a(0)=k(p_0)$.
The inverted stereographic coordinates form a smooth chart at $p_0$,
so $\widehat a$ is $C^2$ near $0$. Moreover, as $p_0$ is a Morse
minimum, $\nabla\widehat a(0)=0$ and $\nabla^2\widehat a(0)$
is positive definite.
Taylor's theorem, applied to $\widehat a$ also to its first derivative, gives
\[
\widehat a(w)=k(p_0)+\frac12\nabla^2\widehat a(0)[w,w]+o(|w|^2),
\]
\[
\nabla\widehat a(w)\cdot w=\nabla^2\widehat a(0)[w,w]+o(|w|^2).
\]
The derivative of $y\mapsto y/|y|^2$ sends the radial vector $y$
to $-y/|y|^2$. Therefore, with $w=y/|y|^2$,
\[
 b(y)=-\nabla\widehat a(w)\cdot w
     =-\nabla^2\widehat a(0)[w,w]+o(|w|^2).
\]
Positive definiteness implies $b(y)\leq-c|y|^{-2}$ for $|y|\geq R_0$,
with some $c>0$ and $R_0>0$. For $R>R_0$, radial integration gives
\[
\int_{R_0<|y|<R}b(y)\,dy
\leq-c\int_{R_0<|y|<R}|y|^{-2}\,dy
\longrightarrow-\infty\qquad \text{as }R\to+\infty,
\]
since $m\geq3$. The integral of $b$ over $B_{R_0}(0)$ is finite,
so one can choose $R_a>R_0$ such that
$b<0$ outside $B_{R_a}(0)\subset\R^m$ and
$\int_{B_{R_a}(0)}b<0$.

For $\lambda\neq0$, by the change of variables in \eqref{leading}, we obtain
\begin{equation}\label{radial-gradient}
 \nabla\widetilde T_k(\lambda,\xi)\cdot(\lambda,\xi)
 =\int_{\R^m} b(y)
       \frac{2^m|\lambda|^m}{(\lambda^2+|y-\xi|^2)^m}\,dy.
\end{equation}
For $\varrho=(\lambda,\xi)$ with $\lambda\neq0$, denote the
maximum and minimum of the function $\frac{2^m|\lambda|^m}{(\lambda^2+|y-\xi|^2)^m}$ in \eqref{radial-gradient}
on the fixed ball $\overline B_{R_a}(0)$ by
$M(\varrho)$ and $m(\varrho)$, respectively. Then it is easy to see that
\[
 \frac{M(\varrho)}{m(\varrho)}\longrightarrow1
 \quad\text{as }\lambda^2+|\xi|^2\longrightarrow\infty,
 \quad\text{uniformly for }\lambda\neq0.
\]
Define $b^+=\max\{b,0\}$ and $b^-=\max\{-b,0\}$, so that
$b=b^+-b^-$. The strict inequality $\int_{B_{R_a}(0)}b<0$
allows us to choose $\delta>0$ with
\[
 (1+\delta)\int_{B_{R_a}(0)}b^+
       -\int_{B_{R_a}(0)}b^-<0.
\]
Hence, when $\lambda^2+|\xi|^2$ is sufficiently large, the integral in \eqref{radial-gradient}
over $B_{R_a}(0)$ is bounded from above by
\[
 M(\varrho)\int_{B_{R_a}(0)}b^+
 -m(\varrho)\int_{B_{R_a}(0)}b^-
 \leq m(\varrho)\left((1+\delta)\int_{B_{R_a}(0)}b^+
                       -\int_{B_{R_a}(0)}b^-\right)<0.
\]
Meanwhile, the integral over $\R^m\setminus B_{R_a}(0)$ is nonpositive.

For $\lambda=0$, we have
$\nabla\widetilde T_k(0,\xi)\cdot(0,\xi)=\omega_m b(\xi)<0$
when $|\xi|$ is large. Together with the preceding estimate for $\lambda\neq0$, this gives
$\nabla\widetilde T_k(\varrho)\cdot\varrho<0$ on
$\partial B_R(0)\subset\R^{m+1}$ for sufficiently large $R$.
The homotopy
$H_t(\varrho)=(1-t)\nabla\widetilde T_k(\varrho)-t\varrho$,
$0\leq t\leq1$, has no boundary zeros, since
\[
H_t(\varrho)\cdot\varrho
=(1-t)\nabla\widetilde T_k(\varrho)\cdot\varrho-tR^2<0
\qquad\text{on }\partial B_R(0).
\]
Hence homotopy invariance yields
\begin{equation}\label{total-degree}
\deg(\nabla\widetilde T_k,B_R(0),0)
=\deg(-\Id,B_R(0),0)
=(-1)^{m+1}.
\end{equation}

Thanks to the above analysis, we see that all critical points of $\widetilde T_k$ lie in $B_R(0)$. The finitely many critical points on $\{\lambda=0\}$ are isolated by \eqref{boundary-hessian};
hence the critical set in $\{\lambda>0\}$ is compact and keeps a distance from $\{\lambda=0\}$. Thus, one can find a smooth relatively compact
open set $\Omega\Subset\Pcal=(0,\infty)\times\R^m$ with no critical point of $T_k$ on its boundary.
Since $\widetilde T_k$ is even in $\lambda$, its gradient has the same
degree on $\Omega$ and its reflection across $\{\lambda=0\}$. Under stereographic projection, $\Crit(a)$ corresponds to
$\Crit(k)\setminus\{p_0\}$, with Morse indices preserved.
Applying the Morse identity to $k$ on $S^m$ and using
$\ind(k,p_0)=0$, we obtain
\[
\sum_{\xi\in\Crit(a)}(-1)^{\ind(a,\xi)}
=\chi(S^m)-1=(-1)^m.
\]
Since the Laplacian signs agree at corresponding critical
points of $a$ and $k$, we deduce that
\begin{align*}
	&\sum_{\xi\in\Crit(a)}
	\operatorname{sgn}(\Delta_{\R^m}a(\xi))(-1)^{\ind(a,\xi)}\\
	&\qquad=\sum_{\xi\in\Crit(a)}(-1)^{\ind(a,\xi)}
	-2\sum_{\substack{\xi\in\Crit(a)\\\Delta_{\R^m}a(\xi)<0}}
	(-1)^{\ind(a,\xi)}\\
	&\qquad=(-1)^m-2S(k).
\end{align*}
And therefore, by additivity, \eqref{boundary-index} and \eqref{total-degree}, we have
\begin{align*}
	(-1)^{m+1}
	&=2\deg(\nabla T_k,\Omega,0)
	+\sum_{\xi\in\Crit(a)}
	\operatorname{sgn}(\Delta_{\R^m}a(\xi))(-1)^{\ind(a,\xi)}\\
	&=2\deg(\nabla T_k,\Omega,0)+(-1)^m-2S(k).
\end{align*}
Solving this identity for the degree gives
$\deg(\nabla T_k,\Omega,0)=S(k)-(-1)^m$.
It is nonzero by assumption \textup{(K2)}, which proves \eqref{degree-formula}.
\end{proof}

\begin{remark}
The preceding degree calculation corresponds to
\cite[Lemmas 3.2-3.3 and Corollary 3.4]{IX}. The first-order
perturbation restricted to the standard spinorial Yamabe family
in that paper is
\[
 -\frac{m-1}{2m}\,m^m
 \int_{\R^m}\frac{a(\lambda y+\xi)}{(1+|y|^2)^m}\,dy.
\]
Our corresponding term is $-(C_m/2)T_k$. 
This agreement concerns the first-order finite-dimensional term;
the coupled linearization and its normal inverse require the
separate analysis in Section~\ref{sec:nondegeneracy}.
\end{remark}

\subsection{Generic nondegeneracy}

\begin{lemma}\label{genericity}
There is a residual subset $\Acal\subset C^2(S^m;\R)$ for which
$T_k$ is Morse on $\Pcal$. This property is independent of the
stereographic pole. For any fixed $p_0$, the intersection
$\Acal\cap\Kcal_{p_0}$ is residual in
\[
 \Kcal_{p_0}=\big\{k\in\Kcal:k(p_0)=\min_{S^m}k\big\}.
\]
\end{lemma}

\begin{proof}
Fix $p_0\in S^m$ as a pole and write $\varrho=(\lambda,\xi)$. Equivalent to
\eqref{leading}, let us consider
\begin{equation}\label{kernel-density}
 T_k(\varrho)=\int_{\R^m}a(y)K_{\varrho}(y)\,dy,\qquad
 K_{\varrho}(y):=K_{\lambda,\xi}(y)
 =\left(\frac{2\lambda}{\lambda^2+|y-\xi|^2}\right)^m.
\end{equation}
Every parameter derivative of $K_{\varrho}$ has an integrable uniform
bound on compact subsets of $\Pcal$. Hence $T_k$ is smooth, even for $k\in C^2$.

Fix $\varrho=(\lambda,\xi)\in\Pcal$ and a nonempty open set $U\subset\R^m$.
For $h\in C_c^\infty(U)$, the variation $a\mapsto a+th$ corresponds
to a smooth perturbation of $k$ on $S^m$.
Since $k\mapsto\nabla T_k(\varrho)$ is linear, its derivative with
respect to $k$ is represented by
\begin{equation}\label{gradient-vector}
h\longmapsto\int_U h(y)\nabla_{\varrho}K_{\varrho}(y)\,dy
\in\R^{m+1}.
\end{equation}
We claim that this map is surjective. If its range were a proper
subspace of $\R^{m+1}$, there would be a nonzero vector
$(\lambda_0,\xi_0)\in\R\times\R^m$ orthogonal to every vector in its range.
Thus
\[
\int_U h(y)\big(\lambda_0\partial_\lambda K_{\varrho}(y)
+\xi_0\cdot\nabla_\xi K_{\varrho}(y)\big)\,dy=0
\qquad\text{for all }h\in C_c^\infty(U).
\]
The arbitrariness of $h$ implies
$\lambda_0\partial_\lambda K_{\varrho}+\xi_0\cdot\nabla_\xi K_{\varrho}\equiv0$
on $U$. Dividing by the positive factor
$mK_{\varrho}/(\lambda^2+|y-\xi|^2)$ gives
\[
\frac{\lambda_0}{\lambda}(|y-\xi|^2-\lambda^2)
+2\xi_0\cdot(y-\xi)=0\qquad \forall y\in U.
\]
A polynomial vanishing on a nonempty open set must be identically zero.
Therefore $(\lambda_0,\xi_0)=(0,0)$, a contradiction. This proves the
claimed surjectivity. In particular, for each fixed $\varrho$, one
can choose $m+1$ such variations of $h\in C_c^\infty(U)$ whose images form a basis of
$\R^{m+1}$; by continuity, they continue to span $\R^{m+1}$ for nearby parameters.

Let $B_n(0)=\{\xi\in\R^m:|\xi|<n\}$ and set
$\mathcal V_n=[1/n,n]\times\overline B_n(0)$. Let $\mathcal U_n$ consist of potentials for which
$T_k$ has no degenerate critical point on $\mathcal V_n$. This is open in the $C^2(S^m,\R)$ topology. In fact, the definition of $T_k$ in \eqref{kernel-density} implies continuity of $k\mapsto T_k$ into $C^2(\mathcal V_n)$.
If $k_j\to k$ and $T_{k_j}$ had degenerate critical points
$\varrho_j\in\mathcal V_n$, compactness of $\mathcal V_n$ would
give a subsequence $\varrho_j\to\varrho\in\mathcal V_n$; uniform
convergence of the first two derivatives would make $\varrho$ a
degenerate critical point of $T_k$. This proves openness.

To see the density of $\mathcal U_n$ in $C^2(S^m,\R)$, let us fix $k$ arbitrarily in $C^2(S^m,\R)\setminus \mathcal U_n$. For each $\varrho_0\in\mathcal V_n$, one can
choose $m+1$ smooth variations of $k$ (compactly supported in $\pi_{p_0}^{-1}(U)\subset S^m$, where
$U\subset\R^m$ is the fixed nonempty open set in \eqref{gradient-vector} and $\pi_{p_0}:S^m\setminus\{p_0\}\to\R^m$ stand for the stereographic projection), generating a basis of $\R^{m+1}$ at $\varrho_0$ via the gradient mapping \eqref{gradient-vector}. By continuity, these vectors remain a basis of $\R^{m+1}$ for parameters
in an open neighborhood of $\varrho_0$. Compactness of $\mathcal V_n$
gives a finite subcover. Collect all the corresponding variation
directions of $k$ and denote them by $h_1,\ldots,h_\ell\in C^\infty(S^m;\R)$.
Then the family $\big\{\nabla T_{h_1}(\varrho),\ldots,\nabla T_{h_\ell}(\varrho)\big\}$ span
$\R^{m+1}$ for every $\varrho$ in an open neighborhood
$\mathcal W\subset\Pcal$ of $\mathcal V_n$. Define 
\begin{equation}\label{perturbation-T}
\mathcal F(t,\varrho)
=\nabla T_{\,k+\sum_{j=1}^{\ell}t_jh_j}(\varrho),
\qquad t=(t_1,\dots,t_\ell)\in\R^\ell,\quad\varrho\in\mathcal W.
\end{equation}
Linearity in the potential gives
\[
\nabla_t\mathcal F(t,\varrho)[\tau]
=\sum_{j=1}^{\ell}\tau_j\nabla T_{h_j}(\varrho),
\qquad\tau\in\R^\ell,
\]
which is surjective onto $\R^{m+1}$ by construction. Hence the full differential of
$\mathcal F$ is surjective, which means that $\mathcal F$ is a submersion on
$\R^\ell\times\mathcal W$, in particular on a neighborhood of
$\{0\}\times\mathcal V_n$.

For ease of notations, let us set $k_t=k+\sum_{j=1}^{\ell}t_jh_j$. Then $\mathcal F:\R^\ell\times \mathcal W\to\R^{m+1}$, $(t,\varrho)\mapsto \mathcal F(t,\varrho)$ is smooth. Since $\mathcal F$ is a submersion, its zero set
\[
\mathcal M=\mathcal F^{-1}(0)
=\big\{(t,\varrho)\in\R^\ell\times\mathcal W:
\mathcal F(t,\varrho)=0\big\}
\]
is a smooth manifold of dimension $\ell$. Consider the restricted
projection
\[
\operatorname{pr}:\mathcal M\longrightarrow\R^\ell,
\qquad\operatorname{pr}(t,\varrho)=t.
\]
At $(t,\varrho)\in\mathcal M$, the tangent space is
\[
T_{(t,\varrho)}\mathcal M
=\big\{(\tau,\vartheta):
\nabla_t\mathcal F(t,\varrho)\cdot\tau
+\nabla_\varrho\mathcal F(t,\varrho)\cdot\vartheta=0\big\}.
\]
The differential of the above projection sends $(\tau,\vartheta)$ to $\tau$,
and $\nabla_\varrho\mathcal F(t,\varrho)=\nabla^2T_{k_t}(\varrho)$. Hence
\[
\ker \nabla\operatorname{pr}\big|_{(t,\varrho)}
=\big\{(0,\vartheta):\nabla^2T_{k_t}(\varrho)\vartheta=0\big\}.
\]
Since both $\mathcal M$ and the coefficient space have dimension
$\ell$, $\nabla\operatorname{pr}\big|_{(t,\varrho)}$ is surjective if and only if its kernel
is trivial, i.e., if $\nabla^2T_{k_t}(\varrho)$ is invertible.
Thus $t$ is a regular value of $\operatorname{pr}$ if and only if
every critical point of $T_{k_t}$ in $\mathcal W$ is nondegenerate.

By Sard's theorem, the critical values of $\operatorname{pr}$ have
Lebesgue measure zero. In particular, every neighborhood of $0$
contains a regular value $t$. For such $t$, we have
$k_t\in\mathcal U_n$, since $\mathcal V_n\subset\mathcal W$, and
\[
\|k_t-k\|_{C^2(S^m)}
\leq\sum_{j=1}^{\ell}|t_j|\,\|h_j\|_{C^2(S^m)}
\to0\qquad \text{as } t\to0 .
\]
Thus $\mathcal U_n$ is dense and
$\Acal=\bigcap_n\mathcal U_n$ is residual.

The construction of $\Acal$ above used a fixed stereographic pole $p_0$ (which is not necessarily a minimum point of a function), we now show that $\Acal$ is independent of this choice, allowing
the pole to be taken at a global minimum of a specified function $k$ as in
Lemma~\ref{degree-lemma}. Changes between stereographic charts are
generated by translations, dilations, orthogonal maps and inversion.
These maps preserve the family of densities $K_{\varrho}(y)\,dy$
and induce diffeomorphisms of $\Pcal$. For the first three maps,
this follows directly from the expressions in \eqref{kernel-density}; for inversion,
direct calculation gives
\[
K_{\lambda,\xi}\!\left(\frac{y}{|y|^2}\right)|y|^{-2m}
=K_{\lambda/(\lambda^2+|\xi|^2),\,\xi/(\lambda^2+|\xi|^2)}(y).
\]
Thus changing the pole replaces $T_k$ by its composition with a
diffeomorphism of $\Pcal$. Such a reparametrization preserves
nondegeneracy of critical points, so elements in $\Acal$ does
not depend on the pole.

It remains to prove that $\Acal\cap\Kcal_{p_0}$ is residual in
$\Kcal_{p_0}$ for a prescribed point $p_0$. Since $\Kcal$ is open in $C^2(S^m,\R)$, $\Kcal_{p_0}$ is relatively
open in the closed (hence complete) set
$\big\{k\in C^2(S^m;\R):k(p_0)=\min_{S^m}k\big\}$, and is therefore a Baire space.

Fix $k\in\Kcal_{p_0}$ and choose a nonempty open set $\mathcal O\subset S^m$
such that $k-k(p_0)>\delta$ on $\mathcal O$ for some $\delta>0$. Repeat the construction
leading to \eqref{perturbation-T}, choosing all variation directions
$h_1,\ldots,h_\ell$ to be compactly supported in $\mathcal O$. For sufficiently small $t$,
we have $\|k_t-k\|_{C^2}<\delta/2$. Since the perturbation vanishes
in $S^m\setminus\mathcal O$, $k_t(p_0)=k(p_0)=\min_{S^m}k_t$.
By openness of $\Kcal$, we also have $k_t\in\Kcal$ for small $t$,
and hence $k_t\in\Kcal_{p_0}$. Applying the preceding Sard argument
with these localized variations shows that
$\mathcal U_n\cap\Kcal_{p_0}$ is dense in $\Kcal_{p_0}$.
Moreover, it is relatively open because $\mathcal U_n$ is open. Finally,
since $\Kcal_{p_0}$ is a Baire space, the Baire category theorem
shows that
\[
\Acal\cap\Kcal_{p_0}
=\bigcap_{n=1}^{\infty}(\mathcal U_n\cap\Kcal_{p_0})
\]
is dense in $\Kcal_{p_0}$. This intersection is residual by
construction, proving the last assertion.
\end{proof}

\subsection{Existence, orbit counts, and separated potentials}

We first record two consequences common to all the constructions.
On a compact parameter domain, the scalar component and spinor
norm of a pair $(v_{\lambda,\xi},\psi_{\lambda,\xi,\ga})\in\Zcal$ have a uniform positive lower
bound. The normal corrections tend to zero in $X$, and hence in
$C^0\times C^0$. Therefore every solution constructed on such a
block satisfies \eqref{positive-conclusion} for small $|\eps|$.
Shrinking $\eps_0$ also ensures $1+\eps k\geq1/2$.
Sobolev embedding and elliptic Schauder estimates give
$v\in C^{2,\alpha}$ and $\psi\in C^{1,\alpha}$ for a suitable
$\alpha>0$, so the solutions are classical.

Testing the two equations of \eqref{system} by $v$ and $\psi$
respectively gives the exact identity
\begin{equation}\label{exact-energy}
 2(m-2)E_\eps(v,\psi)
 =\int_{S^m}(1+\eps k)v^{2/(m-2)}|\psi|^2\dd.
\end{equation}
At $\eps=0$, \eqref{normalization} and conformal invariance
give the constant value of $E_0|_{\Zcal}$ in \eqref{ed-reduction}. Equations
\eqref{ed-reduction} and \eqref{exact-energy} imply
\eqref{energy-asymptotic} for every constructed solution.

\begin{proof}[Proof of the existence assertion in Theorem~\ref{existence}]
Choose $\Omega$ from Lemma~\ref{degree-lemma} and perform the
normal reduction on a neighborhood of
$\overline\Omega\times S^{2N-1}$. The circle action on the unit spinor sphere is free, as verified
in Section~\ref{subsec:mapping property}, and $S^{2N-1}/S^1=\CP^{N-1}$ has Euler
characteristic $\chi(S^{2N-1}/S^1)=\chi(\CP^{N-1})=N\neq0$. The reduced functional is constant on
circle orbits. Denote its induced function on the quotient by
$\overline E_\eps^{\mathrm{red}}$. Equations \eqref{ed-reduction}
and \eqref{leading} give
\[
 \frac{E_0|_{\Zcal}-\overline E_\eps^{\mathrm{red}}(\varrho,[\ga])}{\eps}
 =\frac{C_m}{2}T_k(\varrho)+O(|\eps|)\quad\text{in }C^2.
\]
Apply Lemma~\ref{index-lemma} with
$\mathcal T=(C_m/2)T_k$. It gives an interior critical point for each
sufficiently small nonzero $\eps$. Lifting to the unit spinor sphere
gives a critical point of $E_\eps^{\mathrm{red}}$, since the quotient map is a
submersion and derivatives along the group orbits vanish.
Proposition~\ref{normal-reduction} implies a solution of the
full system. The preceding positivity, regularity, and energy
arguments complete the proof. 
\end{proof}

\begin{proof}[Proof of the orbit counts in Theorem~\ref{existence}]
For $k\in\Acal\cap\Kcal$, Lemmas~\ref{degree-lemma} and
\ref{genericity} imply a nondegenerate critical point $\varrho_0$
of $T_k$. Apply Corollary~\ref{abstract-multiple} at this point,
with $\mathcal T=(C_m/2)T_k$. For the circle action the quotient is
$\CP^{N-1}$, while for the quaternionic action it is
$\HP^{N/2-1}$. Their integral cohomology rings are
\[
 H^*(\CP^{N-1};\mathbb Z)=\mathbb Z[\alpha]/(\alpha^N),
 \quad |\alpha|=2,
 \qquad
 H^*(\HP^{N/2-1};\mathbb Z)=\mathbb Z[\beta]/(\beta^{N/2}),
 \quad |\beta|=4.
\]
The cup-length bound for category gives
$\Cat(\CP^{N-1})\geq N$ and
$\Cat(\HP^{N/2-1})\geq N/2$. The corollary counts distinct
solution orbits and proves the required bounds. The common
arguments above give positivity and the energy estimate.
\end{proof}

\begin{remark}[Orbit distinction and critical values]\label{orbit-distinction}
The preceding proof produces different points of $\CP^{N-1}$ or
$\HP^{N/2-1}$, and the normal-coordinate uniqueness in
Corollary~\ref{abstract-multiple} makes their solution orbits distinct.
These are inequivalent
solutions under the specified internal symmetries, not necessarily
inequivalent metrics under all diffeomorphisms preserving the data.
In particular, in dimension three the bound $N/2=1$ gives only one
quaternionic orbit, not several quaternionically inequivalent solutions. Furthermore, the energy asymptotics do not by themselves distinguish these orbits. Indeed, near a fixed base critical point $\varrho_0$, we have
$\nabla T_k(\varrho_0)=0$ and
$\varrho_\eps(\ga)-\varrho_0=O(|\eps|)$.
Proposition~\ref{normal-reduction} therefore gives, uniformly in $\ga$,
\[
E_\eps^{\mathrm{red}}(\varrho_\eps(\ga),\ga)
=E_0|_{\Zcal}-\eps\frac{C_m}{2}T_k(\varrho_0)
+\eps^2R_0(\varrho_0,\ga)+o(\eps^2).
\]
Here $R_0(\varrho_0,\ga)=\lim_{\eps\to0}R_\eps(\varrho_0,\ga)$,
where $R_\eps$ is defined in \eqref{reduced-expansion}.
The constant and first-order terms are independent of $\ga$,
and $R_0(\varrho_0,\ga)$ may also take the same value at distinct
critical points on the quotient. Thus neither different energies
nor different concentration parameters follow without additional
information. In Theorem~\ref{separated}, different translated
parameter neighborhoods do supply a separate spatial distinction
between the corresponding groups of solution orbits. Under its additional
amplitude condition, those groups also have disjoint energy ranges (as shown below), and this gives at least $l$ energy values.
\end{remark}

\begin{proof}[Proof of Theorem~\ref{separated}]
Put $a_i=k_i\circ\pi_{p_0}^{-1}$ and $T_i=T_{k_i}$.
Linearity gives the exact identity for \eqref{sum-potential},
\begin{equation}\label{sum-leading}
 T_k(\lambda,\xi)=\sum_{j=1}^lT_j(\lambda,\xi-z_j).
\end{equation}
Choose a nondegenerate critical point $\varrho_i$ of $T_i$ and fix
a small neighborhood $U_i\Subset\Pcal$ of $\varrho_i$. Use
$(\lambda,\eta)$ as coordinates in $U_i$, with the corresponding
global parameters be given by $(\lambda,\xi)=(\lambda,z_i+\eta)$.
Keeping $U_i$ fixed, let the centers $z_1,\ldots,z_l$ vary so that
$\min_{1\leq r<s\leq l}|z_r-z_s|\to\infty$. Then, for each $j\neq i$,
\begin{equation}\label{separation-limit}
	\left\|T_j(\cdot,\cdot+z_i-z_j)-\omega_m k_j(p_0)
	\right\|_{C^2(\overline U_i)}\to0.
\end{equation}
Indeed, by \eqref{leading},
\[
T_j(\lambda,\eta+z_i-z_j)
=\int_{\R^m}a_j(\lambda y+\eta+z_i-z_j)\rho(y)\,dy.
\]
On bounded sets of $y$, the argument of $a_j$ tends to infinity
uniformly for $(\lambda,\eta)\in\overline U_i$.
Since $a_j,\nabla a_j,\nabla^2a_j$ tend respectively to
$k_j(p_0),0,0$, the integrand and its parameter derivatives up to
second order converge uniformly to
$k_j(p_0)\rho(y)$ and zero, respectively. Their absolute values
are bounded by $C(1+|y|^2)\rho(y)$, independently of the centers,
because $a_j$ and its derivatives up to second order are bounded.
This bound is integrable for $m\geq3$. Uniform control of the tails,
together with the convergence on bounded $y$ sets and
$\int_{\R^m}\rho=\omega_m$, proves \eqref{separation-limit}.

Write $z_*=(z_1,\ldots,z_l)$. For sufficiently large separation,
\eqref{sum-leading} and \eqref{separation-limit} show that
$T_k(\lambda,z_i+\eta)-\omega_m\sum_{j\neq i}k_j(p_0)$ is a small
$C^2$ perturbation of $T_i(\lambda,\eta)$ on $U_i$.
Let $\varrho_i(z_*)$ denote the resulting nondegenerate critical
point of $T_k$ in the translated neighborhood
$\{(\lambda,z_i+\eta):(\lambda,\eta)\in U_i\}$.
These neighborhoods have disjoint closures for sufficiently large
separation of $z_*$. And in what follows, let us fix such a configuration $z_*$.

For each $i$, consider the compact subset 
\[
\Zcal_i(z_*)=\{z_{\lambda,\xi,\ga}:
(\lambda,\xi)=\varrho_i(z_*),\ \ga\in S^{2N-1}\}\subset\Zcal.
\]
Corollary~\ref{abstract-multiple}, applied with
$\mathcal T=(C_m/2)T_k$ at $\varrho_i(z_*)$, gives at least $N$
distinct $S^1$ solution orbits near $\Zcal_i(z)$, and at least
$N/2$ distinct $S^3$ solution orbits when $m\equiv2,3,4\pmod8$.
This construction uses only the local nondegeneracy of
$\varrho_i(z_*)$; the superposition $k$ defined in \eqref{sum-potential}
need not belong to $\Kcal$ or to $\Acal$.

The compact sets $\Zcal_i(z_*)$ are pairwise disjoint and hence have
positive distance in $X$. Both group actions preserve each set,
since they fix the base parameters $\lambda$ and $\xi$. Moreover, solutions constructed
in different translated neighborhoods cannot belong to the same
orbit: equivariance and the uniform uniqueness of normal coordinates
on compact sets, proved in Corollary~\ref{abstract-multiple}, would
force their base parameters
to coincide, contrary to the disjointness of these neighborhoods.
For each $i$, choose a threshold $\eps_i(z_*)>0$ small enough that
the local construction is valid, its base parameters remain in
the corresponding neighborhood, and this normal-coordinate
uniqueness applies. Then
\[
\eps_0(z_*)=\min_{1\leq i\leq l}\eps_i(z_*)>0
\]
provides a common range $0<|\eps|<\eps_0(z)$ for all $l$
constructions. Their orbit counts therefore add, giving at least
$lN$ distinct $S^1$ orbits and, in the indicated dimensions,
$lN/2$ distinct $S^3$ orbits.

\medskip

For the energy assertion, we first prove a sufficient critical-value
condition that applies to potentials of different shapes.
Write $d=\min_{i\neq j}|z_i-z_j|$ and put
\[
b_i=T_{k_i}(\varrho_i)-\omega_m k_i(p_0).
\]
The $C^2$ convergence in \eqref{separation-limit} and persistence give
\[
T_k(\varrho_i(z_*))
=\omega_m\sum_{j=1}^l k_j(p_0)+b_i+o(1) \quad \text{as } d\to\infty.
\]
If the $b_i$ are pairwise distinct, set
$\delta=\min_{i\neq j}|b_i-b_j|>0$ and choose $d$ sufficiently large
that
\[
|T_k(\varrho_i(z_*))-T_k(\varrho_j(z_*))|\geq\delta/2
\qquad(i\neq j).
\]
Fix such a translation configuration $z_*$, and then choose
$|\eps|$ sufficiently small, depending on $z_*$. It is clear that the parameter branch obtained from Corollary~\ref{abstract-multiple} satisfies
\[
\varrho_{\eps,i}(\ga)
=\varrho_i(z_*)+O_{z_*}(|\eps|),
\qquad \ga\in S^{2N-1},
\]
uniformly in $\ga$. In fact, this follows from the $C^1$ dependence of the
branch on $(\eps,\ga)$ and compactness of the unit spinor sphere.

Let $u_{\eps,i}=(v_{\eps,i},\psi_{\eps,i})\in X$ be a solution of \eqref{system}
obtained by the local multiplicity construction in the $i$-th
region corresponds to some $\ga$. By \eqref{reduced-expansion}, its energy is
\[
\begin{aligned}
	E_\eps(u_{\eps,i})
	&=E_\eps^{\mathrm{red}}(\varrho_{\eps,i}(\ga),\ga)\\
	&=E_0|_{\Zcal}
	-\eps\frac{C_m}{2}T_k(\varrho_{\eps,i}(\ga))
	+\eps^2R_\eps(\varrho_{\eps,i}(\ga),\ga).
\end{aligned}
\]
Since $\nabla T_k(\varrho_i(z_*))=0$ and the Hessian of $T_k$ is
bounded on a fixed compact neighborhood of $\varrho_i(z_*)$,
Taylor's theorem gives
\[
T_k(\varrho_{\eps,i}(\ga))
=T_k(\varrho_i(z_*))+O_{z_*}(\eps^2),
\]
uniformly in $\ga$. Its contribution to the energy error is thus
$O_{z_*}(|\eps|^3)$. Proposition~\ref{normal-reduction} also bounds
$R_\eps$ uniformly on the corresponding compact parameter domain.
Consequently,
\[
E_\eps(u_{\eps,i})
=E_0|_{\Zcal}
-\eps\frac{C_m}{2}T_k(\varrho_i(z_*))
+O_{z_*}(\eps^2).
\]
These bounds hold uniformly over the entire spinorial parameter
space, and hence over all solutions selected by the local
multiplicity construction. Since there are only finitely many
regions, the remainder constants can also be chosen uniformly in
$i$. Of course, they may depend on the fixed configuration $z_*$.

For any solutions selected from distinct regions $i\neq j$, the
 critical-value gap therefore give
\[
|E_\eps(u_{\eps,i})-E_\eps(u_{\eps,j})|
\geq\frac{C_m\delta}{4}|\eps|+O_{z_*}(\eps^2)
\geq\frac{C_m\delta}{8}|\eps|>0
\]
after reducing the admissible range of $|\eps|$. This proves that the energy ranges of the different regions are disjoint, and in particular gives at
least $l$ distinct energy values.

To verify this condition under the amplitude assumption in the theorem,
choose a nondegenerate critical point $\varrho_*$ of $T_h$ and set
$\varrho_i=\varrho_*$ for every $i$. Since $T_{k_i}=t_iT_h$,
these points are nondegenerate and
\[
b_i=t_i\big(T_h(\varrho_*)-\omega_m h(p_0)\big).
\]
The factor in parentheses is strictly positive since $p_0$ is the minimum of $h$. Thus the $b_i$ are pairwise
distinct. Positive multiplication preserves the Morse and degree
conditions defining $\Kcal$, as well as the nondegeneracy of $T_h$.
This proves the stated construction with different energies.
\end{proof}

\section{Noncompactness for a fixed coefficient}\label{noncompactness}

Berti and Malchiodi~\cite{BM01} constructed metrics close
to the round metric for which the Yamabe equation admits
unbounded sequences of positive solutions. Their use of
infinitely many localized perturbations and a uniform
reduction provides a related construction in the scalar
setting. For the present coupled system, we keep the round
metric fixed and follow the translated-potential strategy
of \cite[Section~5]{ISX}.

To construct a noncompact solution sequence for a single
potential $k$, we establish estimates that are uniform after
conformal recentering. We use a compactly supported Euclidean
potential, so that the resulting coefficient is flat to
second order at the concentration point and the interactions
can be controlled relative to the amplitudes of the
translated potentials. A sharper estimate in the natural
energy space controls the nonlinear remainder. Both $k$
and $\eps$ remain fixed along the solution sequence.

\subsection{A compactly supported potential}

For a smooth compactly supported real function $K$ on $\R^m$, write
\begin{equation}\label{compact-seed-term}
 T_K(\lambda,\xi)=\int_{\R^m}K(x)
 \left(\frac{2\lambda}{\lambda^2+|x-\xi|^2}\right)^m\,dx.
\end{equation}
This is the same first-order term as in \eqref{leading}, after a
change of variables.

\begin{lemma}\label{annular-seed}
For every $m\geq3$, there is a nonnegative radial $K\in C_c^\infty(\R^m)$ with $\displaystyle \int K=1$, supported in a thin annulus about the unit sphere, for which $T_K$ has a nondegenerate critical point $\varrho_*=(\lambda_*,0)$ of Morse index one. In particular, there is a small smooth bounded neighborhood $\Omega\Subset\Pcal$ containing no other critical point of $T_K$, with
\[
 \deg(\nabla T_K,\Omega,0)=-1,\qquad
 d_\Omega:=\min_{\partial\Omega}|\nabla T_K|>0.
\]
\end{lemma}

\begin{proof}
Let $\mu$ be the normalized Riemannian volume measure on the unit
sphere $S^{m-1}\subset\R^m$ endowed with its standard round metric,
that is,
\[
d\mu
=\frac{d\mathrm{vol}_{g_{S^{m-1}}}}
{\mathrm{Vol}(S^{m-1},g_{S^{m-1}})}=\frac1{\omega_{m-1}}d\mathrm{vol}_{g_{S^{m-1}}}.
\]
Define $T_\mu:\Pcal=(0,\infty)\times\R^m\to\R$ by
\[
T_\mu(\lambda,\xi)=\int_{S^{m-1}}\left(\frac{2\lambda}{\lambda^2+|x-\xi|^2}\right)^m\,d\mu=\frac1{\omega_{m-1}}\int_{S^{m-1}}\left(\frac{2\lambda}{\lambda^2+|x-\xi|^2}\right)^m\,d\mathrm{vol}_{g_{S^{m-1}}},
\]
where we have embedded $S^{m-1}$ canonically into $\R^m$ and used the coordinates in $\R^m$ for the variable $x\in S^{m-1}$. Then, at $\xi=0$, we have
\[
 T_\mu(\lambda,0)=\left(\frac{2\lambda}{\lambda^2+1}\right)^m.
\]
Taking derivative with respect to $\xi$, in an arbitrary direction $\zeta$, we find that
\[
\nabla_\xi T_\mu(\lambda,\xi)\cdot\zeta
=2m(2\lambda)^m\int_{S^{m-1}}
\frac{(x-\xi)\cdot\zeta}{(\lambda^2+|x-\xi|^2)^{m+1}}\,d\mu.
\]
And at $\xi=0$, this becomes
\[
\nabla_\xi T_\mu(\lambda,0)\cdot\zeta
=\frac{2m(2\lambda)^m}{(\lambda^2+1)^{m+1}}
\int_{S^{m-1}}x\cdot\zeta\,d\mu=0.
\]
This identity holds for every $\lambda>0$, and hence $\partial_\lambda\nabla_\xi T_\mu(\lambda,0)=(0,\dots,0)\in\R^m$ as well. Furthermore, we can deduce that
\[
 T_\mu(1,0)=1,\qquad
 \nabla^2T_\mu(1,0)=\operatorname{diag}(-m,1,\ldots,1).
\]
This gives us a nondegenerate critical point of $T_\mu$ of Morse index one.

Choose $\delta_n\searrow0$, and nonnegative
smooth radial functions $K_n\in C_c^\infty(\R^m)$ such that
\[
\int_{\R^m}K_n(x)\,dx=1,\qquad
\operatorname{supp}K_n
\subset\{x\in\R^m:1-\delta_n<|x|<1+\delta_n\}.
\]
For every bounded continuous function $f:\R^m\to\R$, radial
symmetry of $K_n$ and normalization of the measure $\mu$ give
\[
\int_{\operatorname{supp}K_n}
f\!\left(\frac{x}{|x|}\right)K_n(x)\,dx
=\int_{S^{m-1}}f(x)\,d\mu.
\]
Consequently,
\begin{equation}\label{T-Kn-converge}
\left|\int_{\R^m}f(x)K_n(x)\,dx
	-\int_{S^{m-1}}f(x)\,d\mu\right|\leq
	\sup_{\substack{x\in S^{m-1}\\[0.2em]|r-1|\leq\delta_n}}
	|f(rx)-f(x)|\to0, \quad \text{as } n\to\infty
\end{equation}
by uniform continuity of $f$ on the fixed compact annulus
$\{x:1/2\leq|x|\leq3/2\}$. This proves weak convergence of the
measures $K_n(x)\,dx$ to $d\mu$.

To preserve the nondegenerate critical point at $(1,0)$ of $T_\mu$, we also
need convergence of the first and second parameter derivatives of $T_{K_n}$.
Fix a compact parameter neighborhood of $(1,0)$ contained in
$\Pcal$, and apply the preceding estimate to the kernel $\big(\frac{2\lambda}{\lambda^2+|x-\xi|^2}\big)^m$
in \eqref{compact-seed-term} and to its derivatives with respect
to $(\lambda,\xi)$ up to second order. Since $K_n$ is independent
of these parameters, differentiation under the integral acts only
on the kernel. The kernel and these derivatives are jointly
continuous in $(\lambda,\xi,x)$, and hence uniformly continuous
on the product of the chosen parameter neighborhood with
$\{x:1/2\leq|x|\leq3/2\}$. As $|rx-x|\leq\delta_n$,
similar estimates of \eqref{T-Kn-converge} holds uniformly in $(\lambda,\xi)$
for the kernel and each of these derivatives. Thus $T_{K_n}$ and
its parameter derivatives up to second order converge uniformly
to those of $T_\mu$.

We now use topology degree to obtain a critical point of $T_{K_n}$. Since $\partial_\lambda T_\mu(1,0)=0$ and
$\partial_{\lambda\lambda}T_\mu(1,0)=-m$, choose $0<\eta<1$
sufficiently small so that the $\lambda$-derivative has its only zero
at $\lambda=1$ in $[1-\eta,1+\eta]$, with
\[
\partial_\lambda T_\mu(1-\eta,0)>0,\qquad
\partial_\lambda T_\mu(1+\eta,0)<0.
\]
For all sufficiently large $n$, uniform convergence of $T_{K_n}$ preserves
these endpoint signs. The straight-line homotopy between the
$\lambda$-derivatives of $T_{K_n}$ and $T_\mu$ therefore has no zero
at either endpoint, so we obtain
\[
	\deg\big(\partial_\lambda T_{K_n}(\cdot,0),
	(1-\eta,1+\eta),0\big)=\deg\big(\partial_\lambda T_\mu(\cdot,0),
	(1-\eta,1+\eta),0\big)=-1.
\]
Hence there is $\lambda_*\in(1-\eta,1+\eta)$ with
$\partial_\lambda T_{K_n}(\lambda_*,0)=0$. Radial symmetry of $K_n$ gives
$\nabla_\xi T_{K_n}(\lambda_*,0)=0$, so $(\lambda_*,0)$ is a
critical point of $T_{K_n}$. By choosing $\eta$ sufficiently small and then
$n$ sufficiently large, continuity of $\nabla^2T_\mu$ and the $C^2$
convergence of $T_{K_n}$ ensure that its Hessian has one negative and $m$
positive eigenvalues. Fix such an $n$ and set $K=K_n$.
Choose a sufficiently small smooth bounded neighborhood $\Omega$
of this nondegenerate critical point, containing no other critical
point. The local degree of $\nabla T_K$ is the sign of the Hessian
determinant, hence $-1$, and its gradient is nonzero on
$\partial\Omega$. This proves the lemma.
\end{proof}

\begin{remark}
	Both Lemma~\ref{degree-lemma} and Lemma~\ref{annular-seed}
	establish the existence of critical points of the first-order
	term in the reduced problem. The former obtains existence from
	a nonzero degree, whereas the latter directly constructs a
	potential $K$ for which $T_K$ has a nondegenerate critical point.
	The choice of $K$ here does not rely on conditions
	(K1)-(K2). More precisely, its smooth extension
	to $S^m$, given by $K\circ\pi_{p_0}$, does not
	belong to $\Kcal$. For the local multiplicity argument in Section~\ref{sec:existence}, nondegeneracy
	is ensured by restricting the potential to a residual subset
	of $C^2(S^m,\R)$ provided by Lemma~\ref{genericity}. Here the
	required local nondegeneracy is already guaranteed by the
	construction. Thus no additional generic perturbation is needed
	before applying the local results of Section~\ref{sec:perturbation}, and the chosen
	$K$ can be used directly in the construction below.
\end{remark}

In the sequel, we fix $K$ and $\Omega$ as in the Lemma \ref{annular-seed}, and let $e\in\R^m$ be a unit vector, say for instance $e=(1,0,\dots,0)$. Let us consider the concrete sequences
\begin{equation}\label{fixed-sequences}
 z_i=Ri^2e,\qquad a_i=i^{-10},\qquad i\geq1,
\end{equation}
where $R$ will be chosen large. Define one potential on $S^m$ by
\begin{equation}\label{fixed-potential}
 k(x)=\begin{cases}
 \displaystyle\sum_{i=1}^\infty a_iK(\pi_{p_0}(x)-z_i),&x\neq p_0,\\
 0,&x=p_0.
 \end{cases}
\end{equation}

\begin{lemma}\label{fixed-regularity}
The series \eqref{fixed-potential} converges in $C^2(S^m,\R)$.
It defines a nonnegative nonconstant potential with
\[
 k(p_0)=0,\qquad \nabla_g k(p_0)=0,\qquad \nabla_g^2k(p_0)=0.
\]
\end{lemma}

\begin{proof}
	For each $i$, the function $K(\pi_{p_0}(\cdot)-z_i)$ is supported
	in the compact set
	\[
	\pi_{p_0}^{-1}\big(z_i+\operatorname{supp}K\big)
	\subset S^m\setminus\{p_0\}.
	\]
	It therefore extends by zero to a smooth function on $S^m$ and
	vanishes on a neighborhood of $p_0$. This neighborhood may depend
	on $i$.
	
	We next estimate its $C^2$ norm uniformly in $i$. Fix $L>0$ with
	$\operatorname{supp}K\subset B_L(0)$. Near $p_0$, use the inverted
	stereographic coordinates $w=y/|y|^2$, where $y=\pi_{p_0}(x)$;
	in these coordinates the translated potential is
	$K(w/|w|^2-z_i)$. Wherever this function or one of its derivatives
	is nonzero, we have
	\[
	\left|\frac{w}{|w|^2}-z_i\right|\leq L \quad \text{and} \quad
     |w|^{-1}\leq |z_i|+L.
	\]
	The first and second derivatives of $w\mapsto w/|w|^2$ are bounded
	by $C|w|^{-2}$ and $C|w|^{-3}$, respectively. Since the derivatives
	of $K$ are bounded, the chain rule gives
	\[
	\left|\nabla_w\left[K\!\left(\frac{w}{|w|^2}-z_i\right)\right]\right|\leq C|w|^{-2},
	\]
	\[
		\left|\nabla_w^2\left[K\!\left(\frac{w}{|w|^2}-z_i\right)\right]\right|
		\leq C\big(|w|^{-4}+|w|^{-3}\big).
	\]
	Thus the coordinate derivatives through order two are bounded by
	$C(1+|z_i|)^4$, where the constant $C$ may
	vary from line to line. On the compact complement of a fixed neighborhood
	of $p_0$, the derivatives of $\pi_{p_0}$ are bounded. Hence, we obtain
	\[
	\|K(\pi_{p_0}(\cdot)-z_i)\|_{C^2(S^m)}
	\leq C(1+|z_i|)^4,
	\]
	where $C$ is independent of $i$. By \eqref{fixed-sequences}, we find that
	\[
		\sum_{i=1}^{\infty}a_i
		\|K(\pi_{p_0}(\cdot)-z_i)\|_{C^2(S^m)}
		\leq C\sum_{i=1}^{\infty}i^{-10}(1+Ri^2)^4
		\leq C(1+R)^4\sum_{i=1}^{\infty}i^{-2}<\infty.
	\]
	The series therefore converges in $C^2(S^m)$. Its value, differential,
	and covariant Hessian at $p_0$ can be evaluated term by term.
	Each is zero for every summand, since that summand vanishes near
	$p_0$. Hence $k(p_0)=0$, $\nabla_gk(p_0)=0$, and
	$\nabla_g^2k(p_0)=0$.
\end{proof}

\subsection{Interactions relative to the peak amplitudes}

For $s>1$ set
\[
 S_i(s)=\sum_{j\neq i}a_j|z_i-z_j|^{-s}.
\]
The following estimate will be used with different exponents:
\begin{equation}\label{sequence-estimate}
 S_i(s)\leq CR^{-s}\big(i^{-2s}+i^{-10-s}\big).
\end{equation}
Here $C$ may depend on $s$, but is independent of $i$ and $R$.
To prove this, recall that $a_j=j^{-10}$ and
\[
|z_i-z_j|=R|i^2-j^2|=R|i-j|(i+j).
\]
We split the sum into the three ranges $j\leq i/2$,
$i/2<j<2i$ with $j\neq i$, and $j\geq2i$.

If $j\leq i/2$, then $i^2-j^2\geq3i^2/4$, and hence
\[
\sum_{1\leq j\leq i/2}a_j|z_i-z_j|^{-s}
\leq CR^{-s}i^{-2s}\sum_{j=1}^{\infty}j^{-10}
\leq CR^{-s}i^{-2s}.
\]
If $i/2<j<2i$ and $j\neq i$, then
$a_j\leq2^{10}i^{-10}$ and $|z_i-z_j|\geq Ri|i-j|$.
Each positive integer value of $|i-j|$ occurs at most twice, so
\[
\sum_{\substack{i/2<j<2i\\j\neq i}}a_j|z_i-z_j|^{-s}
\leq CR^{-s}i^{-10-s}\sum_{\ell=1}^{\infty}\ell^{-s}
\leq CR^{-s}i^{-10-s},
\]
where the last series converges because $s>1$.

Finally, if $j\geq2i$, then $j^2-i^2\geq3j^2/4$, giving
\[
\begin{aligned}
	\sum_{j\geq2i}a_j|z_i-z_j|^{-s}&\leq CR^{-s}\sum_{j\geq2i}j^{-10-2s}
	\leq CR^{-s}\int_{2i-1}^{\infty}t^{-10-2s}\,dt\\
	&\leq CR^{-s}i^{-9-2s}
	\leq CR^{-s}i^{-2s}.
\end{aligned}
\]
Adding these three estimates proves \eqref{sequence-estimate}
for every $i\geq1$.

\medskip

To proceed, let
\begin{equation}\label{b-i}
 b_i(\pi_{p_0}^{-1}(y))
 =\sum_j a_jK(y+z_i-z_j),\qquad b_i(p_0)=0.
\end{equation}
This is the coefficient $k$ in coordinates recentered at $z_i$;
equivalently, $b_i=k\circ\Theta_{-z_i}$, where
\[
\Theta_z(x)=
\begin{cases}
	\pi_{p_0}^{-1}\big(\pi_{p_0}(x)-z\big),
	&x\in S^m\setminus\{p_0\},\\
	p_0,&x=p_0,
\end{cases}
\qquad z\in\R^m.
\]
Its corresponding first-order term in our variational framework  is
\[
 T_{b_i}(\lambda,\xi)
 =a_iT_K(\lambda,\xi)
   +\sum_{j\neq i}a_jT_K(\lambda,\xi+z_i-z_j).
\]
Recall that $\Omega\Subset\Pcal$ is the neighborhood of
$(\lambda_*,0)$ fixed in Lemma~\ref{annular-seed}. On
$\overline\Omega$, the scale $\lambda$ is bounded above and away
from zero, and $\xi$ is bounded. Since $K$ has compact support
and $|z_i-z_j|\geq3R$ for $i\neq j$, we may choose $R$ sufficiently
large that
\[
|x-\xi-z_i+z_j|
\geq |z_i-z_j|-|x|-|\xi|
\geq\frac12|z_i-z_j|
\]
for all $x\in\operatorname{supp}K$ and
$(\lambda,\xi)\in\overline\Omega$.

For every multi-index $\nu$ with $|\nu|\leq2$, differentiation under the integral thus gives
\[
\begin{aligned}
	&\big|\partial_{\lambda,\xi}^{\nu}T_K(\lambda,\xi+z_i-z_j)\big|\\
	&\qquad\leq C\int_{\operatorname{supp}K}
	\frac{|K(x)|}{(\lambda^2+|x-\xi-z_i+z_j|^2)^m}\,dx\\
	&\qquad\leq C|z_i-z_j|^{-2m}\int_{\R^m}|K(x)|\,dx
	=C|z_i-z_j|^{-2m}.
\end{aligned}
\]
The constant $C$ depends on the fixed $K$, $\Omega$, and $m$,
but is independent of $i$, $j$, and $R$.
Consequently
\begin{equation}\label{relative-seed-estimate}
 \left\|a_i^{-1}T_{b_i}-T_K\right\|_{C^2(\overline\Omega)}
 \leq C\frac{S_i(2m)}{a_i}
 \leq CR^{-2m}\big(i^{10-4m}+i^{-2m}\big).
\end{equation}
For $m\geq3$ this is uniformly small when $R$ is large, and tends
to zero as $i\to\infty$.

To estimate the nonlinear correction define
\begin{equation}\label{weighted-force}
 \eta_i=a_i+S_i(m+1).
\end{equation}
Inequality \eqref{sequence-estimate} implies
\begin{equation}\label{relative-quadratic}
 \sup_i\frac{\eta_i^2}{a_i}<\infty,\qquad
 \frac{\eta_i^2}{a_i}\to0 \quad \text{as } i\to\infty.
\end{equation}
Indeed the ratio is bounded by
$C\big(i^{-10}+R^{-2(m+1)}i^{6-4m}\big)$, up to smaller terms.
We consider the ratio $\eta_i^2/a_i$ because the normal correction
is $O(|\eps|\eta_i)$ in the natural energy norm, while the resulting
reduced-energy remainder is $O(\eps^2\eta_i^2)$. Since the leading
perturbation term is of order $|\eps|a_i$, the relative size of
the remainder is controlled by $|\eps|\eta_i^2/a_i$.

\begin{remark}
	The specific choice of the exponents in \eqref{fixed-sequences} are only for
	convenience. For each fixed $m\geq3$, the same construction works with
	\[
	a_i=i^{\alpha_1},\qquad z_i=Ri^{\alpha_2}e,
	\qquad \alpha_2>1,\qquad
	-2m\alpha_2<\alpha_1<-4\alpha_2-1.
	\]
	These are sufficient conditions for the estimates used above.
	Indeed, $\alpha_1+4\alpha_2<-1$ gives
	$\sum_i a_i(1+|z_i|)^4<\infty$, which ensures convergence of
	\eqref{fixed-potential} in $C^2(S^m,\R)$ and the vanishing of its
	derivatives through order two at $p_0$. The same
	decomposition used to prove \eqref{sequence-estimate} gives
	\[
	S_i(s)\leq CR^{-s}
	\big(i^{-\alpha_2s}
	+i^{\alpha_1-s(\alpha_2-1)}\big),
	\qquad s>1.
	\]
	Consequently, for $R\geq1$,
	\[
	\frac{S_i(2m)}{a_i}\leq CR^{-2m}
	\big(i^{-\alpha_1-2m\alpha_2}
	+i^{-2m(\alpha_2-1)}\big),
	\]
	\[
		\frac{\eta_i^2}{a_i}
		\leq C\big(i^{\alpha_1}
		+R^{-2(m+1)}
		i^{-\alpha_1-2\alpha_2(m+1)}\big).
	\]
	All powers of $i$ on the right are negative under the stated
	conditions. Thus both ratios are uniformly bounded and tend to
	zero as $i\to\infty$, while the first can be made uniformly small
	by increasing $R$. The uniform reduction and the energy asymptotics
	therefore follow in the same way, yielding a noncompact solution
	sequence with infinitely many distinct energy values as we shall see later. The choice
	$(\alpha_1,\alpha_2)=(-10,2)$ satisfies these conditions for every
	$m\geq3$ and was made only to simplify the estimates.
\end{remark}

\subsection{A uniform reduction after conformal recentering}

Use the fixed compact block
$\overline\Omega\times S^{2N-1}\subset \Pcal\times S^{2N-1}$, and write its corresponding elements
as $z_{\varrho,\ga}\in\Zcal$. Set $\Zcal_c:=\big\{z_{\varrho,\ga}:\varrho\in\overline\Omega,
\ \ga\in S^{2N-1}\big\}$. In this subsection, we shall also use the following natural
energy space:
\[
\mathscr H=H^1(S^m;\R)\times H^{1/2}(\Sigma S^m).
\]
The auxiliary space $H$ of Section~\ref{sec:perturbation} remains the real $L^2$ pairing
space. Throughout this subsection, the $L^p$ norms of $b_i$ (see its definition in \eqref{b-i}) are computed with respect to $d\mathrm{vol}_g$ on $(S^m,g)$, corresponding to
the measure $(2/(1+|y|^2))^m\,dy$ in stereographic coordinates.
The pointwise convergence $b_i\to0$ and uniform boundedness then
yield $\|b_i\|_{L^p}\to0$ by dominated convergence. The finite-rank projection $P_z$ extends boundedly to
$\mathscr H$ and $\mathscr H^*$ by its smooth tangent-frame formula.
On the dual space its action is the transpose of the projection on
$\mathscr H$ under the duality pairing. 
Let $G_b:X_+\to\R$ denote \eqref{perturbation} with $k$ replaced by
$b$, and set
\[
E_{\eps,b}(v,\psi)=E_0(v,\psi)-\eps G_b(v,\psi),
\qquad (v,\psi)\in X_+.
\]
Write $F_{\eps,b}:X_+\to Y$ for the gradient map representing its
first variation as in \eqref{gradient-pairing}; it is obtained from
\eqref{euler-operator} by replacing $k$ with $b$.

\begin{lemma}\label{uniform-fixed-reduction}
There exist $\eps_0>0$ and a neighborhood of $\Zcal_c$ in
$\Zcal$, both independent of $i$, on which the normal corrections
$w_{\eps,i}(z)\in(I-P_z)X$ are defined for $|\eps|<\eps_0$ and satisfy
\[
(I-P_z)F_{\eps,b_i}(z+w_{\eps,i}(z))=0,
\qquad w_{0,i}(z)=0.
\]
Uniformly for $z\in\Zcal_c$, they satisfy
\begin{equation}\label{two-norm-correction}
	\|w_{\eps,i}(z)\|_X\leq C|\eps|\|b_i\|_{L^p},\qquad
	\|w_{\eps,i}(z)\|_{\mathscr H}\leq C|\eps|\eta_i,
\end{equation}
where $\eta_i$ is given in \eqref{weighted-force}.
The scalar component and spinor norm of $z+w_{\eps,i}(z)$ have
uniform positive lower bounds. Define the reduced functional by
\[
E_{\eps,b_i}^{\mathrm{red}}(z)
:=E_{\eps,b_i}(z+w_{\eps,i}(z)).
\]
Expressed in the parameters $(\varrho,\ga)$, it satisfies
\begin{equation}\label{fixed-remainder}
	\left\|
	E_{\eps,b_i}^{\mathrm{red}}(z_{\varrho,\ga})-E_0|_{\Zcal}
	+\eps\frac{C_m}{2}T_{b_i}(\varrho)
	\right\|_{C^1(\overline\Omega\times S^{2N-1})}
	\leq C\eps^2\eta_i^2.
\end{equation}
The constants are independent of $i$ and $|\eps|<\eps_0$.
For each fixed such $\eps$, $\|w_{\eps,i}(z)\|_X\to0$ as
$i\to\infty$, uniformly for $z\in\Zcal_c$.
\end{lemma}

The proof of Lemma \ref{uniform-fixed-reduction} is somewhat lengthy. On a first reading, the reader may proceed directly to the next subsection and return later to verify the details.

\begin{proof}
For each fixed $i$, Proposition~\ref{normal-reduction} applies to
$E_{\eps,b_i}$. However, the perturbation interval and the remainder
bounds obtained in this way may depend on $i$. To construct a sequence
of solutions for a single fixed nonzero $\eps$, we need a common
perturbation interval and estimates for the remainder in the reduced
energy relative to $|\eps|a_i$, uniformly in $i$. We establish these
estimates for the coefficient family $b_i$, using $X$ to construct
the normal corrections and $\mathscr H$ to estimate their contribution
to the reduced energy.

The scalar components and spinor norms of the fields $z\in\Zcal_c$
have uniform positive lower bounds. By compactness and Sobolev
embedding, we can choose a fixed bounded open neighborhood
$U\subset X_+$ of $\Zcal_c$ on which both quantities remain
uniformly bounded away from zero. We carry out the reduction on
a common neighborhood of $\Zcal_c$ in $\Zcal$ contained in $U$, say for instance $U\cap\Zcal$,
and state the estimates on $\Zcal_c$. 
	
The recentered coefficients $b_i$ form a bounded family in
$L^\infty(S^m)$, since
\[
\sup_i\|b_i\|_\infty
\leq\|K\|_\infty\sum_{j=1}^{\infty}a_j<\infty.
\]
The coefficient $b_i$ enters $F_{\eps,b_i}$ only through
multiplication. Hence field differentiation does not differentiate
$b_i$, and the preceding $L^\infty$ bound gives uniform bounds
for the field derivatives on $U$. In particular,
writing $A_z=F_0'(z)$, we have, for $w\in(I-P_z)X$ with
$\|w\|_X\leq r$,
\[
	\|(I-P_z)F_{\eps,b_i}(z)\|_Y\leq C|\eps|\|b_i\|_{L^p},
\]
\[
	\|F_{\eps,b_i}'(z+w)-A_z\|_{\mathcal L(X,Y)}\leq C(r+|\eps|),
\]
uniformly in $i$ and $z\in\Zcal_c$. The first bound
uses $F_0(z)=0$ and the uniform bounds on the standard fields.
The second follows from the uniform second field derivatives of
$F_0$ and the fact that the perturbation is linear in $\eps b_i$.

For $z\in\Zcal_c$, the normal inverses
$A_z^{-1}:(I-P_z)Y\to(I-P_z)X$ are uniformly bounded.
For $w\in(I-P_z)X$ with $\|w\|_X\leq r$, define
\[
\mathcal T_{\eps,i,z}(w)
:=w-A_z^{-1}(I-P_z)F_{\eps,b_i}(z+w).
\]
Its fixed points are exactly the solutions of the normal equation.
Since $A_z^{-1}(I-P_z)A_z$ is the identity on $(I-P_z)X$, the
preceding derivative estimate gives
\begin{equation}\label{Lipschitz}
\|\mathcal T_{\eps,i,z}(w_1)
-\mathcal T_{\eps,i,z}(w_2)\|_X
\leq C(r+|\eps|)\|w_1-w_2\|_X.
\end{equation}
At $w=0$, we have explicitly
\[
\mathcal T_{\eps,i,z}(0)
=-A_z^{-1}(I-P_z)F_{\eps,b_i}(z),\qquad
\|\mathcal T_{\eps,i,z}(0)\|_X
\leq C|\eps|\|b_i\|_{L^p}.
\]
Both bounds hold with a common constant $C$ independent of $i$
and $z\in\Zcal_c$.

Choose $r>0$ sufficiently small that $Cr\leq1/4$ and
$z+w\in U$ whenever $z\in\Zcal_c$ and $\|w\|_X\leq r$.
This is possible because $\Zcal_c$ is compact and $U$ is open.
Next, using $\sup_i\|b_i\|_{L^p}<\infty$, choose $\eps_0>0$ so that
\[
C\eps_0\leq\frac14,\qquad
C\eps_0\sup_i\|b_i\|_{L^p}\leq\frac r2.
\]
For $|\eps|<\eps_0$, the Lipschitz constant in \eqref{Lipschitz} is therefore at most
$1/2$, and $\|\mathcal T_{\eps,i,z}(0)\|_X\leq r/2$.
In particular, for $\|w\|_X\leq r$,
\[
\|\mathcal T_{\eps,i,z}(w)\|_X\leq \|\mathcal T_{\eps,i,z}(w)-\mathcal T_{\eps,i,z}(0)\|_X+\|\mathcal T_{\eps,i,z}(0)\|_X
\leq\frac12\|w\|_X+\frac r2\leq r.
\]
Thus $\mathcal T_{\eps,i,z}$ maps the closed normal $X$-ball of
radius $r$ into itself and is a contraction there. As in
Proposition~\ref{normal-reduction}, it has a unique fixed point
$w_{\eps,i}(z)$, with the same choices of $r$ and $\eps_0$ for
every $i$ and $z\in\Zcal_c$. The fixed-point identity gives
\[
\|w_{\eps,i}(z)\|_X
\leq \frac12\|w_{\eps,i}(z)\|_X+C|\eps|\|b_i\|_{L^p},
\]
and hence
\[
\|w_{\eps,i}(z)\|_X\leq C|\eps|\|b_i\|_{L^p}.
\]
This proves the first estimate in \eqref{two-norm-correction}.
Since $z+w_{\eps,i}(z)\in U$, the uniform positive scalar and
spinor norm lower bounds are also preserved.

To proceed, let us mention that the operator $A_z$ is also an isomorphism from the normal subspace
of $\mathscr H$ to its dual normal subspace. Its principal part is
$\operatorname{diag}(L_g,D_g)$, an isomorphism from $\mathscr H$
to $\mathscr H^*$, while its lower-order terms are compact by the compact
embeddings into $L^2$. Every $h\in\mathscr H$ satisfying $A_zh=0$ weakly is smooth by
elliptic regularity, so Proposition~\ref{kernel} gives
$\ker_{\mathscr H} A_z=T_z\Zcal$. Symmetry implies that $\Ran A_z$ is contained
in the orthocomplement of $T_z\Zcal$ in $\mathscr H^*$. Since $A_z$
has index zero, its range has codimension $\dim T_z\Zcal$ and
therefore coincides with this orthocomplement. Restricting $A_z$ to
$(I-P_z)\mathscr H$ removes the kernel and yields the claimed
normal isomorphism. In local trivializations, the normal inverses
depend continuously on $z$ and are therefore uniformly bounded
for $z\in\Zcal_c$.

We next estimate $G_{b_i}'$ and its derivative in tangent
directions. Let $z\in\Zcal_c$ and $u=z+w=(v,\psi)\in U$. Following
the conformal identification of Section~\ref{sec:nondegeneracy}, define the Euclidean fields
\[
\left\{
\begin{aligned}
	\widetilde v(y)
	&=f_0(y)^{(m-2)/2}v(\pi_{p_0}^{-1}(y)),\\[0.3em]
	\widetilde\psi(y)
	&=f_0(y)^{(m-1)/2}\beta^{-1}
	\big(\psi(\pi_{p_0}^{-1}(y))\big),
\end{aligned}\right.
\qquad f_0(y)=\frac{2}{1+|y|^2}.
\]
Here $\beta$ is the conformal identification of the spinor bundles. The uniform bounds on $U$ give
\[
c f_0(y)^{(m-2)/2}\leq\widetilde v(y)
\leq C f_0(y)^{(m-2)/2},\qquad
|\widetilde\psi(y)|\leq C f_0(y)^{(m-1)/2}.
\]
For $\vartheta\in T_z\Zcal$ with $\|\vartheta\|_X\leq1$, the
same upper bounds hold for the scalar and spinor components of
its conformal transform, uniformly for $z\in\Zcal_c$. The two Euclidean components representing $G_b'(u)$ are
\[
\frac{b\circ\pi_{p_0}^{-1}}{m-2}
\widetilde v^{(4-m)/(m-2)}|\widetilde\psi|^2,
\qquad
(b\circ\pi_{p_0}^{-1})\widetilde v^{2/(m-2)}\widetilde\psi.
\]
They are bounded respectively by
$C|b\circ\pi_{p_0}^{-1}|f_0^{(m+2)/2}$ and
$C|b\circ\pi_{p_0}^{-1}|f_0^{(m+1)/2}$.
For $m\geq5$, the lower bound on $\widetilde v$ controls the
negative powers in these estimates. Differentiating the two
expressions in the direction $\vartheta$ preserves the same
powers of $f_0$.

For $b=b_i$, fix $j\neq i$ and consider the term
$a_jK(y+z_i-z_j)$. On its support, we may have
\[
|y|\geq\frac12|z_i-z_j|,
\qquad f_0(y)\leq C|z_i-z_j|^{-2},
\]
provided $R$ is sufficiently large in the definition of $z_i$ in \eqref{fixed-sequences}. The scalar component of $G_{b_i}'(u)$ therefore satisfies
\[
\left|
\frac{a_jK(y+z_i-z_j)}{m-2}
\widetilde v^{(4-m)/(m-2)}|\widetilde\psi|^2
\right|
\le Ca_j|z_i-z_j|^{-(m+2)}.
\]
The spinor component satisfies
\[
	\left|a_jK(y+z_i-z_j)
	\widetilde v^{2/(m-2)}\widetilde\psi\right|
	\leq Ca_j|z_i-z_j|^{-(m+1)}.
\]

To estimate the tangent derivatives, write $\vartheta=(\sigma,\eta)$
and denote its conformal transform by
$(\widetilde\sigma,\widetilde\eta)$. As noted above,
\[
|\widetilde\sigma|\leq C f_0^{(m-2)/2},\qquad
|\widetilde\eta|\leq C f_0^{(m-1)/2}.
\]
Differentiation in the direction $\vartheta$ acts only on the field
variables, so the coefficient $a_jK(y+z_i-z_j)$ remains unchanged.
The product rule gives
\[
\left|\partial_\vartheta
\big(\widetilde v^{(4-m)/(m-2)}|\widetilde\psi|^2\big)\right|\leq C\left(
\widetilde v^{(6-2m)/(m-2)}|\widetilde\sigma|
|\widetilde\psi|^2
+\widetilde v^{(4-m)/(m-2)}
|\widetilde\psi||\widetilde\eta|\right)
\leq C f_0^{(m+2)/2}
\]
and
\[
\left|\partial_\vartheta
	\big(\widetilde v^{2/(m-2)}\widetilde\psi\big)\right|\leq C\left(
	\widetilde v^{(4-m)/(m-2)}|\widetilde\sigma|
	|\widetilde\psi|
	+\widetilde v^{2/(m-2)}|\widetilde\eta|\right)
	\leq C f_0^{(m+1)/2}.
\]
Thus the scalar and spinor components of the $j$-th contribution in $G_{b_i}'$,
as well as their derivatives in the direction $\vartheta$, are
bounded respectively by a constant multiplication of $a_j|z_i-z_j|^{-(m+2)}$ and
$a_j|z_i-z_j|^{-(m+1)}$.

To pass to the dual norm, recall the Sobolev embeddings
\[
H^1(S^m)\hookrightarrow L^{2m/(m-2)}(S^m),\qquad
H^{1/2}(\Sigma S^m)\hookrightarrow
L^{2m/(m-1)}(\Sigma S^m).
\]
By duality, these give the continuous embeddings
\[
L^{2m/(m+2)}(S^m)\hookrightarrow H^{-1}(S^m),\qquad
L^{2m/(m+1)}(\Sigma S^m)\hookrightarrow
H^{-1/2}(\Sigma S^m).
\]
Thus a pair $(f_1,f_2)$ in these Lebesgue spaces represents an
element of $\mathscr H^*$, and H\"older's inequality yields
\[
\|(f_1,f_2)\|_{\mathscr H^*}
\leq C\left(
\|f_1\|_{L^{2m/(m+2)}}+
\|f_2\|_{L^{2m/(m+1)}}\right).
\]
All translated copies of $K$ have the same Euclidean support
volume, the term $a_iK$
contributes at most $Ca_i$, and hence we obtain
\[
\begin{aligned}
	\|G_{b_i}'(u)\|_{\mathscr H^*}
	&\leq C\left[a_i+\sum_{j\neq i}a_j\left(
	|z_i-z_j|^{-(m+2)}+|z_i-z_j|^{-(m+1)}\right)\right]\\
	&\leq C\big(a_i+S_i(m+1)\big)=C\eta_i.
\end{aligned}
\]
Here the second inequality uses $|z_i-z_j|\geq1$, which holds
for sufficiently large $R$. Since the tangent derivatives satisfy
the same pointwise estimates, the same argument applies to
$G_{b_i}''(u)[\vartheta,\cdot]$. Combining the two bounds gives
\begin{equation}\label{force-dual-bound}
	\|G_{b_i}'(u)\|_{\mathscr H^*}
	+\|G_{b_i}''(u)[\vartheta,\cdot]\|_{\mathscr H^*}
	\leq C\eta_i.
\end{equation}
The constant $C$ depends only on the fixed data and is independent
of $i$, $z\in\Zcal_c$, $u=z+w\in U$, and
$\vartheta\in T_z\Zcal$ with $\|\vartheta\|_X\leq1$.
By linearity in $\vartheta$, the second term is bounded by
$C\eta_i\|\vartheta\|_X$ for an arbitrary
$\vartheta\in T_z\Zcal$.

Since $E_0'(z)=0$ for all $z\in\Zcal_c$, Taylor's formula gives
\[
E_0'(z+w)=A_zw+\mathcal N_z(w),
\]
where
\[
\mathcal N_z(w)
:=\int_0^1\big(E_0''(z+tw)-E_0''(z)\big)[w,\cdot]\,dt.
\]
Here the second-order derivatives are interpreted as bounded bilinear forms on
$\mathscr H$, so the integral defines an element of $\mathscr H^*$.
On the small $X$ neighborhood, we shall have
\begin{equation}\label{Ncal-w}
\|\mathcal N_z(w)\|_{\mathscr H^*}
\leq C\|w\|_X\|w\|_{\mathscr H}.
\end{equation}
Next, for $w=w_{\eps,i}(z)\in(I-P_z)X$, the normal equation can be
written in $\mathscr H^*$ as
\[
(I-P_z)\big(E_0'(z+w)-\eps G_{b_i}'(z+w)\big)=0.
\]
Since
$(I-P_z)A_zw=A_zw$, substitution of the expansion above gives
\[
A_zw=(I-P_z)\big(\eps G_{b_i}'(z+w)-\mathcal N_z(w)\big).
\]
The normal inverse
\[
A_z^{-1}:(I-P_z)\mathscr H^*\longrightarrow(I-P_z)\mathscr H
\]
satisfies $\|A_z^{-1}f\|_{\mathscr H}\leq C\|f\|_{\mathscr H^*}$,
uniformly for $z\in\Zcal_c$. Applying this inverse and using the
uniform boundedness of $I-P_z$, \eqref{force-dual-bound}, and the
estimate for $\mathcal N_z(w)$, we obtain
\[
	\|w\|_{\mathscr H}
	\leq C\left(
	|\eps|\|G_{b_i}'(z+w)\|_{\mathscr H^*}
	+\|\mathcal N_z(w)\|_{\mathscr H^*}\right)
	\leq C|\eps|\eta_i+C\|w\|_X\|w\|_{\mathscr H}.
\]
Absorb the last term using the uniform small $X$-norm estimate proves
the second estimate in \eqref{two-norm-correction}. Moreover, the normal equation gives
\[
E_{\eps,b_i}'(z+w)=P_zE_{\eps,b_i}'(z+w).
\]
Since $P_zA_zw=0$, the expansion above yields
\[
E_{\eps,b_i}'(z+w)
=P_z\big(\mathcal N_z(w)-\eps G_{b_i}'(z+w)\big).
\]
Using the uniform boundedness of $P_z$ on $\mathscr H^*$,
\eqref{force-dual-bound}, and the estimate \eqref{Ncal-w},
we obtain
\[
\begin{aligned}
	\|E_{\eps,b_i}'(z+w)\|_{\mathscr H^*}
	&\leq C\left(
	\|\mathcal N_z(w)\|_{\mathscr H^*}
	+|\eps|\|G_{b_i}'(z+w)\|_{\mathscr H^*}\right)\\[0.3em]
	&\leq C\|w\|_X\|w\|_{\mathscr H}+C|\eps|\eta_i\leq C|\eps|\eta_i.
\end{aligned}
\]

We now verify \eqref{fixed-remainder}. Recalling that
$G_{b_i}(z_{\varrho,\ga})=\frac{C_m}{2}T_{b_i}(\varrho)$. By the definition of the reduced functional,
\[
\begin{aligned}
	&E_{\eps,b_i}^{\mathrm{red}}(z)-E_0(z)+\eps G_{b_i}(z)\\
	&\qquad=E_0(z+w)-E_0(z)
	-\eps\big(G_{b_i}(z+w)-G_{b_i}(z)\big).
\end{aligned}
\]
Since $E_0'(z)=0$ for $z\in\Zcal_c$, the linear term in the Taylor expansion vanishes. Hence
\[
	|E_0(z+w)-E_0(z)|=\left|\int_0^1(1-t)E_0''(z+tw)[w,w]\,dt\right|\leq C\|w\|_{\mathscr H}^2,
\]
due to $E_0''(u)$ extends to a uniformly bounded bilinear form on $\mathscr H$ for $u\in U$.
For the perturbation term, the fundamental theorem of calculus
and \eqref{force-dual-bound} give
\[
\begin{aligned}
	|G_{b_i}(z+w)-G_{b_i}(z)|
	&=\left|\int_0^1G_{b_i}'(z+tw)[w]\,dt\right|\\
	&\leq\int_0^1\|G_{b_i}'(z+tw)\|_{\mathscr H^*}
	\|w\|_{\mathscr H}\,dt\\
	&\leq C\eta_i\|w\|_{\mathscr H}.
\end{aligned}
\]
Combining these estimates and using the second bound in
\eqref{two-norm-correction}, we obtain
\[
	\left|E_{\eps,b_i}^{\mathrm{red}}(z)-E_0(z)+\eps G_{b_i}(z)\right|\leq C\|w\|_{\mathscr H}^2
	+C|\eps|\eta_i\|w\|_{\mathscr H}\leq C\eps^2\eta_i^2.
\]

To estimate the derivative of the remainder, let
$\vartheta\in T_z\Zcal$ with $\|\vartheta\|_X\leq1$.
Since $\vartheta\in\ker A_z$, symmetry gives
$(A_zw)[\vartheta]=0$. Together with $E_0'(z)=0$, Taylor's formula
therefore yields
\[
E_0'(z+w)[\vartheta]
=\int_0^1(1-t)E_0'''(z+tw)[w,w,\vartheta]\,dt.
\]
The remaining terms contain no spatial derivatives of the variations,
and their coefficients are uniformly bounded for $z+tw\in U$, while Sobolev embedding gives
$\|\vartheta\|_{C^0}\leq C$. Hence
\[
|E_0'''(z+tw)[w,w,\vartheta]|
\leq C\|\vartheta\|_{C^0}\|w\|_{L^2}^2
\leq C\|w\|_{\mathscr H}^2,
\]
and consequently
$|E_0'(z+w)[\vartheta]|\leq C\|w\|_{\mathscr H}^2$.

For the perturbation term, the fundamental theorem of calculus
and the symmetry of $G_{b_i}''$ allow us to place $\vartheta$
in the first argument and apply \eqref{force-dual-bound}:
\[
\begin{aligned}
	\big|(G_{b_i}'(z+w)-G_{b_i}'(z))[\vartheta]\big|
	&=\left|\int_0^1G_{b_i}''(z+tw)[\vartheta,w]\,dt\right|\\
	&\leq\int_0^1
	\|G_{b_i}''(z+tw)[\vartheta,\cdot]\|_{\mathscr H^*}
	\|w\|_{\mathscr H}\,dt\\
	&\leq C\eta_i\|w\|_{\mathscr H}.
\end{aligned}
\]

The dependence of $w=w_{\eps,i}(z)$ on $z$ contributes one further term through
the chain rule. Differentiating $P_zw=0$ gives
\[
P_z\,dw[\vartheta]=-(dP_z[\vartheta])w.
\]
Since $E_{\eps,b_i}'(z+w)$ annihilates normal directions, only
this tangential component of $dw[\vartheta]$ contributes. Hence
\[
	\big|E_{\eps,b_i}'(z+w)[dw[\vartheta]]\big|
	=\big|E_{\eps,b_i}'(z+w)[P_z\,dw[\vartheta]]\big|
	\leq C\|E_{\eps,b_i}'(z+w)\|_{\mathscr H^*}
	\|w\|_{\mathscr H}
	\leq C\eps^2\eta_i^2.
\]
Combining these estimates with the chain rule gives a bound of
$C\eps^2\eta_i^2$ for the derivative of
$E_{\eps,b_i}^{\mathrm{red}}-E_0|_{\Zcal}
+\eps G_{b_i}|_{\Zcal}$.
This proves
\eqref{fixed-remainder}.

Finally, for each fixed $|\eps|<\eps_0$, the convergence
$\|b_i\|_{L^p}\to0$ established at the beginning of this subsection
and the first estimate in \eqref{two-norm-correction} imply
\[
\sup_{z\in\Zcal_c}\|w_{\eps,i}(z)\|_X\longrightarrow0.
\]
This completes the proof.
\end{proof}

\begin{remark}
	The purpose of the estimate in $\mathscr H$ is to control the
	remainder relative to the leading perturbation of size $|\eps|a_i$.
	The convergence $\|b_i\|_{L^p}\to0$ alone does not provide this
	comparison. By \eqref{fixed-remainder}, the relative $C^1$ error
	is bounded by
	\[
	C|\eps|\frac{\eta_i^2}{a_i}.
	\]
	In view of \eqref{relative-quadratic}, this bound is uniformly
	small for sufficiently small $|\eps|$ and tends to zero as
	$i\to\infty$ for fixed $\eps$. Together with
	\eqref{relative-seed-estimate}, this allows the same nonzero
	$\eps$ to be used in the degree argument for every $i$.
	For a fixed $\eps$, the energy remainder is also $o(a_i)$, which
	will be used below to distinguish infinitely many critical values.
\end{remark}

\subsection{Existence in every parameter block and concentration}

We prove Theorem~\ref{fixed-noncompactness} for the potential
\eqref{fixed-potential}. The $S^1$-symmetry of system \eqref{system} permits us to work on
$Q=S^{2N-1}/S^1=\CP^{N-1}$, with $\chi(Q)=N>0$.
For $\eps\neq0$ define the normalized reduced functional
\[
 \mathcal F_{\eps,i}
 =\frac{2(E_0|_{\Zcal}-E_{\eps,b_i}^{\mathrm{red}})}{C_m\eps a_i}
 \quad\text{on }\overline\Omega\times Q.
\]
Equations \eqref{relative-seed-estimate} and \eqref{fixed-remainder}
give
\begin{equation}\label{uniform-normalized-energy}
 \|\mathcal F_{\eps,i}-T_K\|_{C^1}
 \leq C\frac{S_i(2m)}{a_i}
       +C|\eps|\frac{\eta_i^2}{a_i}.
\end{equation}
First choose $R$ large, then choose one $\eps_0>0$ small, so that
this bound is less than $d_\Omega/2$ for every $i$ and
$0<|\eps|<\eps_0$, where $d_\Omega$ has been given in Lemma \ref{annular-seed}. For each $i$, the gradient field $\nabla\mathcal F_{\eps,i}$ on $\Omega\times Q$, taken with respect to the product metric,
is homotopic to $(\nabla T_K,0)$ through vector fields that do
not vanish on $\partial\Omega\times Q$. The total index is
\[
 \deg(\nabla T_K,\Omega,0)\chi(Q)=-N\neq0,
\]
by the argument of Lemma~\ref{index-lemma}. Thus
$\mathcal F_{\eps,i}$ has an interior critical point
$(\varrho_{\eps,i},[\ga_{\eps,i}])$ for every $i$.
Lift it to the spinor sphere. Invariance makes derivatives along
the orbit vanish, so Proposition~\ref{normal-reduction} gives a full
solution $\widehat u_{\eps,i}$ of the recentered equation.
Translate its Euclidean fields back by $z_i$ to obtain
$u_{\eps,i}=(v_{\eps,i},\psi_{\eps,i})$ for the single coefficient
$1+\eps k$. Conformal covariance preserves the functional and
pointwise positivity of both components $v_{\eps,i}$ and $|\psi_{\eps,i}|$. 

For a fixed $\eps\neq0$, the right-hand side of
\eqref{uniform-normalized-energy} tends to zero as $i\to\infty$.
Since $\Omega$ contains only the critical point $\varrho_*$ of $T_K$,
criticality implies $\varrho_{\eps,i}\to\varrho_*$.
Also $\|w_{\eps,i}\|_X\to0$. The normalized solutions are therefore
asymptotic, in $X$, to the compact standard family at $\varrho_*$.
If $\widehat U_{\eps,i},\widehat\Phi_{\eps,i}$ denote their Euclidean
representatives, uniformly in $i$ they satisfy the pointwise estimates
\begin{equation}\label{fixed-profile-bound}
 \left\{
 \aligned
 &c f_0(y)^{(m-2)/2}\leq\widehat U_{\eps,i}(y)
 \leq C f_0(y)^{(m-2)/2},\\[0.3em]
 &c f_0(y)^{(m-1)/2}\leq|\widehat\Phi_{\eps,i}(y)|
 \leq C f_0(y)^{(m-1)/2}.
 \endaligned \right.
\end{equation}
The pull-back fields satisfy
\[
 \begin{split}
 v_{\eps,i}(\pi_{p_0}^{-1}(x))
 &=f_0(x)^{-(m-2)/2}\widehat U_{\eps,i}(x-z_i),\\
 |\psi_{\eps,i}(\pi_{p_0}^{-1}(x))|
 &=f_0(x)^{-(m-1)/2}|\widehat\Phi_{\eps,i}(x-z_i)|.
 \end{split}
\]
At $x=z_i$ these give lower bounds
$c(1+|z_i|^2)^{(m-2)/2}$ and
$c(1+|z_i|^2)^{(m-1)/2}$. Both suprema tend to infinity.
For bounded sets of $x$  the fields tend uniformly to zero, which gives
local convergence away from $p_0$.

The energy asymptotics are sharper than a uniform $O(\eps^2)$ bound.
Using \eqref{relative-seed-estimate}, \eqref{relative-quadratic}, and
$\varrho_{\eps,i}\to\varrho_*$ in \eqref{fixed-remainder} gives,
for each fixed nonzero sufficiently small $\eps$,
\begin{equation}\label{fixed-critical-values}
 E_\eps(u_{\eps,i})
 =\frac{C_m\omega_m}{2(m-2)}
 -\eps a_i\frac{C_m}{2}T_K(\varrho_*)+o(a_i)
 \qquad(i\to\infty).
\end{equation}
Since $K\geq0$ and $\displaystyle \int K=1$, we have $T_K(\varrho_*)>0$.
The critical values converge to $\frac{C_m\omega_m}{2(m-2)}$ from below
if $\eps>0$ and from above if $\eps<0$. They consequently include infinitely many distinct values.

Testing the equations as in \eqref{exact-energy} now gives
\[
 \int_{S^m}(1+\eps k)v_{\eps,i}^{2/(m-2)}|\psi_{\eps,i}|^2\,
 d\mathrm{vol}_g\to C_m\omega_m.
\]
To prove the measure convergence in \eqref{fixed-concentration-measure}, let
\[
q_i(y)=(1+\eps b_i(\pi_{p_0}^{-1}(y)))
\widehat U_{\eps,i}(y)^{2/(m-2)}|\widehat\Phi_{\eps,i}(y)|^2.
\]
The conformal factors in
$v_{\eps,i}^{2/(m-2)}|\psi_{\eps,i}|^2\,d\mathrm{vol}_g$
cancel, since $f_0^{-1}f_0^{-(m-1)}f_0^m=1$.
Changing variables $x=y+z_i$ and using
$k(\pi_{p_0}^{-1}(y+z_i))=b_i(\pi_{p_0}^{-1}(y))$, we obtain,
for every $f\in C(S^m,\R)$,
\[
\int_{S^m}f\cdot (1+\eps k)v_{\eps,i}^{2/(m-2)}|\psi_{\eps,i}|^2\,
d\mathrm{vol}_g
=\int_{\R^m}f(\pi_{p_0}^{-1}(y+z_i))q_i(y)\,dy.
\]
Taking $f\equiv1$ and using the total mass limit above gives
\[
\int_{\R^m}q_i(y)\,dy\to C_m\omega_m\qquad \text{as } i\to\infty.
\]
Moreover, \eqref{fixed-profile-bound} and the uniform bounds on the
coefficients imply
\[
0\leq q_i(y)\leq C(1+|y|^2)^{-m}.
\]
In particular, their masses outside a large ball are uniformly small:
\[
\sup_i\int_{|y|>L}q_i(y)\,dy
\leq C\int_{|y|>L}(1+|y|^2)^{-m}\,dy
\to0\qquad\text{as }L\to\infty.
\]
For each fixed $L>0$, the fact that $|z_i|\to\infty$ gives
$\pi_{p_0}^{-1}(y+z_i)\to p_0$ uniformly for $|y|\leq L$.
Continuity of $f$ then implies uniform convergence of $f(\pi_{p_0}^{-1}(y+z_i))$ to $f(p_0)$ for $|y|\leq L$. Splitting the integral into
$|y|\leq L$ and $|y|>L$, we obtain
\[
\int_{S^m}f\cdot (1+\eps k)v_{\eps,i}^{2/(m-2)}|\psi_{\eps,i}|^2\,
d\mathrm{vol}_g
\to C_m\omega_mf(p_0),
\]
which proves the weak measure convergence in \eqref{fixed-concentration-measure}.

Finally, we would like to close the proof by mentioning that the $\mathscr H$-bounds for $(v_{\eps,i},\psi_{\eps,i})$ follows from the facts
\[
\int_{S^m}v_{\eps,i}L_gv_{\eps,i}\dd
=\frac1{m-2}\int_{S^m}(1+\eps k)
v_{\eps,i}^{2/(m-2)}|\psi_{\eps,i}|^2\dd
\to\frac{C_m\omega_m}{m-2}>0
\]
and
\[
\|\psi_{\eps,i}\|_{H^{1/2}}
\leq C\|D_g\psi_{\eps,i}\|_{H^{-1/2}}
\leq C
\|v_{\eps,i}\|_{L^{2m/(m-2)}}^{2/(m-2)}
\|\psi_{\eps,i}\|_{L^{2m/(m-1)}}\leq C,
\]
where the last inequality follows from the $H^1$-boundedness of $v_{\eps,i}$ and the pointwise estimates in \eqref{fixed-profile-bound}. All these complete the proof of Theorem~\ref{fixed-noncompactness}.

\section{A Kazdan-Warner obstruction}\label{sec:obstruction}

It is natural to ask whether \eqref{system} admits a nontrivial
solution with positive scalar component for every
$k\in C^2(S^m,\R)$, without the restrictions imposed above.
We shall exhibit a smooth Morse function for which no such solution
exists for any nonzero $\eps$. The argument
is motivated by the obstruction and the height-function example
of Ammann, Humbert, and Ould Ahmedou~\cite{AHOA}.
For the coupled system, the relevant identity involves the
density $v^{2/(m-2)}|\psi|^2$.

\begin{proposition}\label{coupled-kw}
	Let $k\in C^2(S^m,\R)$ and $\eps\neq0$. If $(v,\psi)$ is a
	classical solution of \eqref{system} with $v>0$, then
	\begin{equation}\label{kw-identity}
		\int_{S^m}V(k)v^{2/(m-2)}|\psi|^2\dd=0
	\end{equation}
	for every conformal Killing vector field $V$ on $(S^m,g)$.
\end{proposition}

\begin{proof}
Let $\{F_t\}$ be the flow generated by $V$, characterized by
\[
\frac{\partial}{\partial t}F_t(x)=V(F_t(x)),
\qquad F_0(x)=x.
\]
Since $V$ is a conformal Killing vector field, each $F_t$
is a conformal diffeomorphism. Thus we may write
\[
F_t^*g=\kappa_t^2g,\qquad \kappa_t>0,\quad\kappa_0\equiv1.
\]
Denote by $\mathcal B_t$ the isometric spinorial pullback
obtained by lifting this flow from the identity. Consider
the conformal transforms
\[
v_t=\kappa_t^{(m-2)/2}(v\circ F_t),
\qquad
\psi_t=\kappa_t^{(m-1)/2}\mathcal B_t\psi.
\]
Since $v$ is strictly positive on the compact sphere, elliptic regularity provides sufficient regularity to
differentiate this curve in the functional. Its scalar component remains positive.
	
The conformal covariance of $L_g$ and $D_g$ gives
\[
E_0(v_t,\psi_t)=E_0(v,\psi).
\]
Moreover, the conformal weights satisfy
\[
v_t^{2/(m-2)}|\psi_t|^2
=
\kappa_t^m
\big(v^{2/(m-2)}|\psi|^2\big)\circ F_t.
\]
As $F_t^*\dd=\kappa_t^m\dd$, a change of variables yields
\[
E_\eps(v_t,\psi_t)
=
E_0(v,\psi)
-\frac{\eps}{2}
\int_{S^m}(k\circ F_{-t})
v^{2/(m-2)}|\psi|^2\dd.
\]
Consider the pair $(v,\psi)$ is a critical point of $E_\eps$ and differentiate at $t=0$, we obtain
\[
0
=
\left.\frac{d}{dt}\right|_{t=0}E_\eps(v_t,\psi_t)
=
\frac{\eps}{2}
\int_{S^m}V(k)v^{2/(m-2)}|\psi|^2\dd.
\]
Since $\eps\neq0$, this proves \eqref{kw-identity}.
\end{proof}

\begin{corollary}\label{height-obstruction}
	Let $h(x)=x_{m+1}$, for $x\in S^m$, where $S^m$ is embedded canonically into $\R^{m+1}$ so that $x$ has coordinates $x=(x_1,\dots,x_{m+1})$.
	For $k=h$ and every $\eps\neq0$, system \eqref{system}
	has no classical solution with $v>0$.
	In particular, this nonexistence conclusion holds for
	$0<|\eps|<1$, when the coupling coefficient
	$1+\eps h$ is everywhere positive.
\end{corollary}

\begin{proof}
	On the unit round sphere,
	\[
	\nabla_g^2h=-hg,\qquad
	\Delta_gh=-mh,\qquad
	|\nabla_gh|^2=1-h^2.
	\]
	Thus $V=\nabla_gh$ is a conformal Killing vector field.
	Applying Proposition~\ref{coupled-kw} gives
	\[
	\int_{S^m}(1-h^2)v^{2/(m-2)}|\psi|^2\dd=0.
	\]
	The integrand is nonnegative, and $1-h^2$ is strictly
	positive away from the two poles. Since $v>0$ and
	$\psi$ is continuous, it follows that $\psi\equiv0$.
	The scalar equation then gives $L_gv=0$. Consequently,
	\[
	0=\int_{S^m}vL_gv\dd
	=\int_{S^m}
	\left(|\nabla_gv|^2+\frac{m(m-2)}4v^2\right)\dd,
	\]
	which contradicts $v>0$.
\end{proof}

The height function also permits a direct comparison with
the hypotheses of our existence theorems. Its only critical
points are
\begin{equation}\label{p+-}
p_+=(0,\dots,0,1),\qquad p_-=(0,\dots,0,-1).
\end{equation}
At these points,
\[
\nabla_g^2h(p_+)=-g,\qquad
\nabla_g^2h(p_-)=g.
\]
Hence $h$ is Morse, with
\[
\ind(h,p_+)=m,\qquad \ind(h,p_-)=0,
\qquad
\Delta_gh(p_+)=-m,\qquad \Delta_gh(p_-)=m.
\]
It satisfies \textup{(K1)}, whereas
\[
S(h)
=
\sum_{\substack{\xi\in\Crit(h)\\\Delta_gh(\xi)<0}}
(-1)^{\ind(h,\xi)}
=
(-1)^m.
\]
Thus \textup{(K2)} fails. Corollary~\ref{height-obstruction}
shows that \textup{(K2)} cannot simply be omitted from a
theorem asserting existence for every potential satisfying
\textup{(K1)}.

In fact, this obstruction persists under sufficiently
small $C^2$ perturbations of the height function.

\begin{proposition}\label{open-obstruction}
	Let $h$ be the height function defined above. There exists
	$\delta>0$ such that, for every $k\in C^2(S^m,\R)$ satisfying
	\[
	\|k-h\|_{C^2(S^m)}<\delta,
	\]
	system \eqref{system} has no classical solution with $v>0$
	for any $\eps\neq0$. Moreover, $\delta$ can be chosen so that
	$1+\eps k>0$ on $S^m$ for all such $k$ whenever $|\eps|<1/2$.
\end{proposition}

\begin{proof}
	For $k$ sufficiently close to $h$ in $C^2$, it has
	exactly two critical points $q_+$ and $q_-$, close to
	$p_+$ and $p_-$ in \eqref{p+-}, respectively.
	
	We claim that there is a conformal diffeomorphism $F$, close to the
	identity, such that
	\[
	F(p_+)=q_+,\qquad F(p_-)=q_-.
	\]
	To see this, let $e_1,\ldots,e_{m+1}$ be the standard basis
	of $\R^{m+1}$. For $j=1,\ldots,m$, consider the vector fields
	\[
	R_j(x)=x_{m+1}e_j-x_je_{m+1},
	\qquad
	G_j(x)=\nabla_gx_j=e_j-x_jx.
	\]
	The fields $R_j$ generate rotations in the
	$(x_j,x_{m+1})$ planes and are Killing vector fields.
	The fields $G_j$ are conformal Killing vector fields,
	since $\nabla_g^2x_j=-x_jg$. At the two poles, $p_+=e_{m+1}$ and $p_-=-e_{m+1}$, 
	\[
	R_j(p_+)=e_j,\qquad R_j(p_-)=-e_j,
	\qquad
	G_j(p_+)=G_j(p_-)=e_j.
	\]
	For $a,b\in\R^m$, let $F_{a,b}$ be the time-one map of the flow generated by the conformal Killing vector field $\sum_{j=1}^m(a_jR_j+b_jG_j)$.
	As $a,b$ vary in $\R^m$, these are conformal diffeomorphisms depending smoothly
	on $(a,b)$, with $F_{0,0}=\mathrm{id}$.
	Consider the map
	\[
	(a,b)\longmapsto
	\big(F_{a,b}(p_+),F_{a,b}(p_-)\big)
	\in S^m\times S^m.
	\]
	Identifying both tangent spaces at the poles with $\R^m$,
	its differential at $(0,0)$ is
	\[
	(a,b)\longmapsto(a+b,-a+b).
	\]
	This is clearly an isomorphism, and the inverse function theorem therefore implies that,
	for every ordered pair $(q_+,q_-)$ sufficiently close
	to $(p_+,p_-)$, there are small parameters $(a,b)$ such that
	\[
	F_{a,b}(p_+)=q_+,\qquad F_{a,b}(p_-)=q_-.
	\]
	We choose $F=F_{a,b}$.
	
	Set $f=k\circ F$. Then $f$ is close to $h$ in $C^2$,
	and both poles are critical points of $f-h$.
	Integrating its Hessian along a shortest geodesic
	from the nearer pole gives
	\[
	|\nabla_g(f-h)(x)|
	\leq
	\|\nabla_g^2(f-h)\|_\infty
	d_g(x,\{p_+,p_-\}).
	\]
	If we set $r=d_g(x,\{p_+,p_-\})$, then $0\leq r\leq\pi/2$
	and $|\nabla_gh(x)|=\sin r$. Since
	$r\leq(\pi/2)\sin r$, we obtain
	\[
	|\nabla_g(f-h)|
	\leq
	\frac{\pi}{2}\|\nabla_g^2(f-h)\|_\infty
	|\nabla_gh|.
	\]
	By reducing $\delta$, we may assume that
	$\|\nabla_g^2(f-h)\|_\infty<2/\pi$. Then it follows immediately that
	\[
	(\nabla_gh)(f)=g(\nabla_gh,\nabla_gf)=|\nabla_gh|^2
	+g(\nabla_gh,\nabla_g(f-h))
	\geq
	\left(1-\frac{\pi}{2}
	\|\nabla_g^2(f-h)\|_\infty\right)|\nabla_gh|^2
	>0
	\]
	away from the poles.
	
	The pushforward $V=F_*(\nabla_gh)$ is a conformal
	Killing vector field and satisfies
	\[
	(V(k))\circ F=(\nabla_gh)(k\circ F)>0
	\]
	away from $p_+$ and $p_-$. Therefore $V(k)$ is
	strictly positive away from $q_+$ and $q_-$.
	Identity \eqref{kw-identity} forces $\psi\equiv0$,
	and the positivity of $L_g$ gives the same
	contradiction as in Corollary~\ref{height-obstruction}.
	
	Finally, reducing $\delta$ so that $\|k\|_\infty<2$
	ensures $1+\eps k>0$ whenever $|\eps|<1/2$. This completes the whole proof.
\end{proof}

\begin{remark}
	The residual set $\Acal$ in Lemma~\ref{genericity}
	describes generic nondegeneracy of the reduced problem.
	For $k\in\Kcal$, the existence of critical points of $T_k$
	follows from the nonzero degree in Lemma~\ref{degree-lemma}.
	The additional restriction $k\in\Acal$ ensures their
	nondegeneracy and yields the stated orbit lower bounds.
	This genericity statement does not extend existence to
	a residual subset of the whole space $C^2(S^m,\R)$.
	Indeed, Proposition~\ref{open-obstruction} provides a
	nonempty open subset of $C^2(S^m,\R)$ consisting of
	potentials for which no classical solution with $v>0$
	exists. Consequently, for each fixed sufficiently small
	nonzero $\eps$, the set of potentials admitting such
	a solution is not dense and therefore cannot be residual
	in $C^2(S^m,\R)$.
\end{remark}

\section*{AI Declaration Statement}
All mathematical ideas, arguments, and results presented in this manuscript are solely due
to the authors, and have been obtained without the use of generative AI. Chat GPT was partly exploited in the revision process, towards notational uniformization and streamlined exposition.

\end{document}